\documentclass{amsart}

\usepackage{amssymb}
\usepackage[all]{xy}

\newtheorem{thm}{Theorem}

\newtheorem{lem}[thm]{Lemma}
\newtheorem{cor}[thm]{Corollary}

\newtheorem{prop}[thm]{Proposition}

\theoremstyle{definition}
\newtheorem{defn}[thm]{Definition}

\newtheorem{say}[thm]{}
\newtheorem{exmp}[thm]{Example}

\newtheorem{rem}[thm]{Remark}          

\newtheorem*{ack}{Acknowledgments}      

\newtheorem{defn-thm}[thm]{Definition--Theorem}  
\newtheorem{defn-lem}[thm]{Definition--Lemma}  

\theoremstyle{remark}

\renewcommand{\c}[0]{{\mathbb C}}  

\renewcommand{\o}[0]{{\mathcal O}} 

\renewcommand{\r}[0]{{\mathbb R}} 

\renewcommand{\a}[0]{{\mathbb A}}

\newcommand{\p}[0]{{\mathbb P}}
\newcommand{\f}[0]{{\mathbb F}}

\newcommand{\qtq}[1]{\quad\mbox{#1}\quad}

\newcommand{\coker}[0]{\operatorname{coker}}

\newcommand{\sing}[0]{\operatorname{Sing}}

\newcommand{\rdown}[1]{\lfloor{#1}\rfloor}

\newcommand{\rest}[0]{\operatorname{rest}}

\newcommand{\tsum}[0]{\textstyle{\sum}}

\def\into{\DOTSB\lhook\joinrel\to}

\def\loccoh#1.#2.#3.#4.{H^{#1}_{#2}(#3,#4)}

\DeclareMathAlphabet{\mathchanc}{OT1}{pzc}%
                                {m}{it}

\newcommand{\PSU}{\mathrm{PSU}}

\usepackage[all]{xy}\xyoption{dvips}

\newcommand{\tprod}[0]{\textstyle{\prod}} 

\begin{document}
\bibliographystyle{amsalpha}


\title[Szemer\'edi--Trotter theorems]{Szemer\'edi--Trotter-type theorems\\
 in dimension 3 }
\author{J\'anos Koll\'ar}
\begin{abstract}
We estimate the number of incidences in a configuration of
$m$ lines and $n$ points in dimension 3. The main term is $mn^{1/3}$
if we work over the real or complex numbers but   $mn^{2/5}$
over finite fields. 
\end{abstract}

\maketitle

Let ${\mathcal L}$ be a set of  lines 
and  ${\mathcal P}$  a set of  points in
some affine or projective space.
The papers \cite{MR729791, MR2047237, el-sh, MR2763049, MR2946447} and 
 \cite{gu-ka1, gu-ka, 2014arXiv1404.2321G}
 point out the importance of
bounding the number of special points of  ${\mathcal L}\cup {\mathcal P}$.

\begin{defn} \label{i.i.defn}
Let ${\mathcal L}$ be a set of  lines in some ambient space. 
There are at least 3 sensible way to count the  {\it number of intersection
points} of ${\mathcal L}$. The smallest of these is to count each
 intersection point with multiplicity 1.
 Our formulas give naturally a larger number
$$
I({\mathcal L}):=\tsum_{p}  \bigl(r(p)-1\bigr),
\eqno{(\ref{i.i.defn}.1)}
$$
where $r(p)$ denotes the number of lines
passing through $p$ and 
the summation is over all points where at least 2 lines meet.
The largest is to count all pairs of lines that intersect, this
corresponds to $\tsum_{p}  \binom{r(p)}{2}$.

In addition, let  ${\mathcal P}$ be  a set of  points. 
 An {\it incidence} is a pair  
$(p\in \ell)$ where $\ell\in {\mathcal L}$ and 
$p\in {\mathcal P}$. The total number of incidences
is denoted by  $I({\mathcal L}, {\mathcal P})$. Thus
$$
I({\mathcal L},{\mathcal P})=\tsum_{p\in{\mathcal P} } r(p).
\eqno{(\ref{i.i.defn}.2)}
$$
\end{defn}

The  Szemer\'edi--Trotter theorem \cite{MR729791} says that
for $m$ lines and $n$ points in $\r^2$ the number of incidences satisfies
$$
I({\mathcal L}, {\mathcal P})\leq 2.5 m^{2/3}n^{2/3}+m+n,
\eqno{(\ref{i.i.defn}.3)}
$$
where the  constants are due to \cite{MR1464571, MR2267545}. 
This implies the same bound in any $\r^d$
since projecting to $\r^2$ can only increase the number of incidences.
Thus it is interesting to look for other bounds for
point/line configurations in $\r^d$ or $\c^d$ that do not hold
for planar  ones. 
Any finite configuration of lines and points 
can be projected to $\r^3$  (resp.\ to  $\c^3$)
without changing the number of incidences, 
hence it is enough to study  lines and points in 3--space.

The bounds we obtain are   not symmetric in $m,n$.
Note, however, that while lines and points have  symmetric roles in $\r^2$,
they do not have symmetric behavior in $\r^3$. 

 \begin{thm} \label{main.thm.2}
Let ${\mathcal L}$ be a set of $m$ distinct lines 
and ${\mathcal P}$  a set of $n$ distinct points
in $\c^3$.
Let $c$ be a constant such that no plane  contains
more than $c\sqrt{m}$ of the lines.
Then the number of incidences satisfies 
$$
I({\mathcal L}, {\mathcal P})
\leq (3.66+0.91c^2) m n^{1/3} + 6.76 n.
\eqno{(\ref{main.thm.2}.1)}
$$
\end{thm}

\begin{exmp} Choose ${\mathcal P}$ to be the integral points in the cube
$[0,r-1]^3$ and ${\mathcal L}$ to be the lines parallel to one of the 
coordinate axes meeting ${\mathcal P}$. Then $m=3r^2$, $n=r^3$,
 $I({\mathcal L}, {\mathcal P})=3r^3$ and each plane contains at most
$2r$ lines.  

We can do slightly better by tilting the above configuration.
This is obtained by first taking the image of  ${\mathcal L}$ in $\r^7$
under the map  $(x,y,z)\mapsto (xyz, xy,yz,zx, x,y,z)$
and then projecting generically to $\r^3$. The incidences are unchanged
but now any plane contains at most 2 of the lines.
We can thus take $c$ arbitrary small and (\ref{main.thm.2}.1) gives that
$3r^3=I({\mathcal L}, {\mathcal P})\leq 17.74 r^3$.
Hence, for this series of examples, 
  (\ref{main.thm.2}.1) is   a factor of $< 6$  
away from an optimal bound.

In this example  $mn^{1/3}$ and $n$ are of the same size and the proof
seems to work naturally in this case. Other cases are discussed in
Paragraph  \ref{campare.say}.
\end{exmp}

A key step of the proof of Theorem \ref{main.thm.2}
does not work over finite fields and we have the following
estimate in general. 

 \begin{thm} \label{main.thm.2.fq}
Let ${\mathcal L}$ be a set of $m$ distinct lines 
and ${\mathcal P}$  a set of $n$ distinct points
in $K^3$ for an arbitrary field $K$.
Let $c$ be a constant such that no plane  contains
more than $c\sqrt{m}$ of the lines.
Then the number of incidences satisfies
$$
I({\mathcal L}, {\mathcal P})
\leq 
2.45mn^{2/5}+2.45n^{6/5} +0.91c^2mn^{1/3}+6.74n.
\eqno{(\ref{main.thm.2.fq}.1)}
$$
\end{thm}

 Example \ref{d4.lines.charp.exmp} gives a line/point configuration over
$\f_{q^2}$ where  
 $m=(q+1)(q^3+1) $,   $n=(q^2+1)(q^3+1) $  
and $I({\mathcal L}, {\mathcal P})=(q+1)(q^2+1)(q^3+1)$. For this
(\ref{main.thm.2.fq}.1) gives an upper bound 
$2.45 q^4\bigl(q^5\bigr)^{2/5}+2.45\bigl(q^5\bigr)^{6/5}$
+ lower terms.
Therefore
$$
q^6\leq I({\mathcal L}, {\mathcal P})
\leq 4.9 q^6+(\mbox{lower terms}).
$$
Hence  (\ref{main.thm.2.fq}.1) is   a factor of $< 5$  
away from an optimal bound.
\medskip

Theorems \ref{main.thm.2}--\ref{main.thm.2.fq} give the following form
 of Bourgain's conjecture,
proved in \cite{gu-ka1} over $\c$. 
 As pointed out in \cite{2013arXiv1311.1479E}, 
the exponent $5/4$ is optimal over finite fields; see also
Example \ref{d4.lines.charp.exmp}.

\begin{cor} \label{main.thm.bourg.2}
Let ${\mathcal L}$ be a set of $m$ distinct lines 
and ${\mathcal P}$  a set of $n$ distinct points
in $K^3$.
Assume  that every line contains at least  $\sqrt{m}$ points
and no plane  contains
more than $\sqrt{m}$  of the lines. 
\begin{enumerate}
\item If $K$ has characteristic 0  then 
$n\geq \tfrac1{50}\cdot m^{3/2}$.
\item If $K$  has positive characteristic and $m\geq 10^4$ then 
$n\geq \tfrac1{20}\cdot m^{5/4}$.
\end{enumerate}
\end{cor}

We get somewhat worse bounds for  $I({\mathcal L})$. Note that
 $I({\mathcal L})$ is the largest  when any 2 lines meet
in distinct points; then we get
$\binom{m}{2}$ intersection points. If this  happens
then all the lines are contained in a plane. One gets a similar quadratic 
growth if all lines are  contained in a quadric surface. 
To avoid these  cases, one should assume that no
plane or quadric contains  too many of the lines. 
The following is a strengthening of  \cite[Thm.2.10]{gu-ka},
which in turn was 
 conjectured by \cite{el-sh}.

\begin{thm} \label{main.thm}
Let ${\mathcal L}$ be a set of $m$ distinct lines in $\c^3$.
Let $c$ be a constant such that no plane (resp.\ no quadric) contains
more than $c\sqrt{m}$ (resp.\ more than $2c\sqrt{m}$) of the lines.
Then the number of  intersection points---with multiplicity as in
(\ref{i.i.defn}.1)---is  
$$
I({\mathcal L})\leq \bigl(29.1+\tfrac{c}2\bigr)\cdot m^{3/2}.
$$
\end{thm}

Note that 
\cite[Thm.2.10]{gu-ka} gives a similar bound where each
intersection point is counted with multiplicity 1. 
On the other hand, \cite[Thm.2.11]{gu-ka} gives a bound,
under some restrictions, where each
intersection point is counted with multiplicity $r(p)^2$;
see Remark \ref{2.11.thm.rem} for details.
For the applications in \cite{gu-ka} the relevant value of $c$ is
$\leq 3.5$. 

\begin{say}[Comparison with previous results]\label{campare.say}
The idea of using algebraic surfaces to attack such problems
is due to \cite{gu-ka1}.  
Our main observation is that the arithmetic genus of a
line configuration provides a  very efficient way to
bound the number of incidences.

It was observed by Ellenberg and Hablicsek  as well as by Guth and Katz
(both unpublished) that, at least over $\r$, estimates similar to
Theorem \ref{main.thm.2} could be deduced from the results of
\cite{gu-ka}, though the resulting constants were never computed. 
(In this area, some proofs lead to quite large constants. 
For example, the coefficient 2.5 in (\ref{i.i.defn}.3)
 first appeared as $\leq 10^{60}$, and the complex version,
 due to \cite{2003math5283T},
still has a  coefficient  $\leq 10^{60}$.)

To compare the various bounds, write  $n=m^t$. For 
Theorem \ref{main.thm.2} the interesting cases are when each line contains at
least 2 points and  through each point there are at least
2 lines. Thus $t\in [\frac12, 2]$ and Theorem \ref{main.thm}
shows that in fact $t\in [\frac12, \frac32]$ is the important range.
The easy planar bound 
$I({\mathcal L}, {\mathcal P})\leq  m^{1/2}n$ 
(see Paragraph \ref{planar.say}) is better than
(\ref{main.thm.2}.1) for  $t\in [\frac12, \frac34]$. 
Thus Theorem \ref{main.thm.2} gives new results when
$t\in [\frac34, \frac32]$ and is sharp at $t=\frac32$. 

For real lines, \cite{gu-ka} yields bounds that are smaller
for $t<\frac32$. I do not know whether these  bounds also hold
over $\c$, but  it is  not clear what should replace the ham sandwich theorems
 used in  \cite{gu-ka} to get a proof over $\c$ or over finite fields.

By contrast, the proof of 
Theorem \ref{main.thm} given in \cite{gu-ka}
also works for complex lines, 
hence over any field of characteristic 0. 
Thus the  new part of Theorem \ref{main.thm} 
is the explicit constant. 
The same applies to 
 Corollary \ref{main.thm.bourg.2}.1.

In all these results, the different behavior in positive
characteristic is restricted to small vales of $p$. 
This was observed in  \cite{2013arXiv1311.1479E}
where a positive
characteristic version of   Corollary \ref{main.thm.bourg.2} is also proved.
 We show in (\ref{39.say}--\ref{main.thm.char.p}) that
 Theorem \ref{main.thm}
holds over a field of characteristic $p$ provided
$p>\sqrt{m}$  and Theorem \ref{main.thm.2} 
 holds provided $p>\sqrt[3]{6n}$; answering a question posed by Dvir
in a conversation.


\cite{MR2946447} proves
 higher dimensional analogs of the  Szemer\'edi--Trotter theorem 
where lines are replaced by  larger linear spaces.

The methods of this paper have been applied in \cite{2014arXiv1408.5791H} to
estimate the number of joints in higher dimensional line arrangements.

\end{say}

\begin{say}[Planar case]\label{planar.say}
For planar configurations we clearly have
$$
\tsum \binom{r(p)}{2}=\binom{m}{2},
\eqno{(\ref{planar.say}.1)}
$$
where we sum over all intersection points, thus
$\tsum \bigl(r(p)-1\bigr)^2< m^2$. 
Using the Cauchy--Schwartz inequality
as in (\ref{pf.of.main.thm.2.5}.2),
 this implies that
$I({\mathcal L}, {\mathcal P})\leq  mn^{1/2}$.
This is quite sharp since $m$ general lines in a plane
have $n=\binom{m}{2}$
intersection points and for these
$I({\mathcal L}, {\mathcal P})\sim \bigl(1/\sqrt{2}\bigr)  mn^{1/2}$.
Working with the dual configuration gives that
$I({\mathcal L}, {\mathcal P})\leq  m^{1/2}n$
and the two together imply that
$$
I({\mathcal L}, {\mathcal P})\leq  m^{3/4}n^{3/4}.
\eqno{(\ref{planar.say}.2)}
$$
This is weaker than (\ref{i.i.defn}.3). Note, however, that
(\ref{planar.say}.2) holds over any field and it is sharp over finite fields.
If ${\mathcal L}$ is the set of all lines
and $ {\mathcal P}$ is the set of all points  over $\f_q$ then
$$
q(q^2+q+1)=I({\mathcal L}, {\mathcal P})\leq  m^{3/4}n^{3/4}
= (q^2+q+1)^{3/4}(q^2+q+1)^{3/4}
\eqno{(\ref{planar.say}.3)}
$$
shows that in (\ref{planar.say}.2)  the exponents  $3/4$ and
the constant factor $1$ are all optimal.
\end{say}

\begin{say}[Outline of the proofs]\label{outline} 
For all the theorems there are
 4  steps, the first two follow \cite{gu-ka1, gu-ka}.

(\ref{outline}.1) By an easy dimension count, all the lines 
and points lie on 
a low degree algebraic surface $S$. In general $S$ is reducible;
it is not hard to deal with the components that contain
infinitely many lines. We recall the needed results in
Section \ref{ruled.sec}.

(\ref{outline}.2) In the remaining cases,  old results of Monge, 
Salmon and Cayley are used to
find another surface of low degree $T$ that contains  all the lines.
Thus  the union of all lines  
$C:=\cup \{\ell:\ell\in {\mathcal L}\}$
is contained in the complete intersection curve $B:=S\cap T$.
Since the references are not easily accessible,
we outline the proofs in Section \ref{sketch.sec}.

(\ref{outline}.3) Sometimes we find a lower degree surface $T$. 
Some of the proofs
work without this step but it improves the bounds 
substantially.

(\ref{outline}.4)
 Although $B$ is a singular algebraic curve, the expected formula bounds
its arithmetic genus, hence also
the arithmetic genus  of $C$. 
The key fact is that while a plane curve of degree $d$ has genus
 $\asymp d^2$, a typical complete intersection   curve of degree $d$
in $\p^3$ has genus
 $\asymp d^{3/2}$.
Finally the set of intersection points
equals the set of singular points of $C$ which in turn is
controlled by the arithmetic genus of $C$.

These steps appear in the cleanest form in the proof of
Theorem \ref{main.thm}; we treat it  in Section \ref{sec.intersections}.
The proofs of Theorems \ref{main.thm.2} and \ref{main.thm.2.fq},
presented in  Sections \ref{sec.incidences.C}--\ref{sec.incidences.Fq},
are slightly more involved.
\end{say}

\begin{ack} I thank Z.~Dvir,  J.~Ellenberg, L.~Guth, N.~Katz, S.~Rams, F.~Russo,
  M.~Sharir, J.~Solymosi  and F.~Voloch for
comments, discussions and references. I am especially grateful to 
 S.~Kleiman 
for calling my attention to several papers on the Gauss map and
Hermitian varieties, 
P.~Yang for telling me about \cite{monge-old}
and C.~Stibitz for simplifying the arguments in Section
\ref{coh.sec}. I also thank the referee for many very helpful suggestions.  

Partial financial support   was provided  by  the NSF under grant numbers 
DMS-0968337 and DMS-1362960.
\end{ack}

\section{Low degree surfaces}

The following elementary lemmas show that any collection of lines
or points 
is contained in a relatively low degree algebraic surface.
Under some extra conditions there are  even 2 such surfaces.
We work in projective 3-space $\p^3$ over an arbitrary field.
For a set of lines  ${\mathcal L}$ let
 $[{\mathcal L}]\subset \p^3$ denote their union. We view
$[{\mathcal L}]$ as a (reducible) algebraic curve in $\p^3$.

\begin{lem} \label{lowder.surf.lem}
 Let ${\mathcal L}$ be $m$ distinct lines in $\p^3$.
\begin{enumerate}
\item There is a surface $S$ of degree $d\leq \sqrt{6m}-2$
that contains  $[{\mathcal L}]$.
\item Let $U\subset \p^3$ be an irreducible surface of degree $g\leq \sqrt{6m}$.
Then there is a surface $T$ of degree $e$ 
that contains  $[{\mathcal L}]$,
does not contain $U$ and $ge\leq 6m$. 
\end{enumerate}
\end{lem}

Proof.  Degree $d$ homogeneous polynomials in 4 variables
form a vector space of  dimension $\binom{d+3}{3}$.
For a surface of degree $d$ it is $d+1$ linear conditions to
contain a line. Thus if $\binom{d+3}{3}>m(d+1)$ then
such a surface $S$ exists, giving (1).

For $e\geq g$ the equations of surfaces that contain $U$ form a vector space 
of dimension $\binom{e-g+3}{3}$. Thus if
$$
\tbinom{e+3}{3}-\tbinom{e-g+3}{3}> m(e+1)
$$
then we find   a surface $T$ of degree $e$ 
that contains all the lines in ${\mathcal L}$ but
does not contain $U$. By expanding we see that
$$
\tbinom{e+3}{3}-\tbinom{e-g+3}{3}>\tfrac16 g(e+1)(e+5),
$$
so we are done if $g(e+5)\geq 6m$
since a vector space  can not be a union of $\leq 2$ lower dimensional
vector subspaces. Finally note that if 
$e=\rdown{\frac{6m}{g}}\geq  \frac{6m}{g}-1$
then $g(e+5)\geq g\bigl(\frac{6m}{g}+4\bigr)=6m+4g$.
\qed

\begin{lem} \label{lowder.surf.pts.lem}
 Let ${\mathcal P}$ be $n$ distinct points in $\p^3$.
\begin{enumerate}
\item There is a surface $S$ of degree $d\leq \sqrt[3]{6n}$
that contains  ${\mathcal P}$.
\item Let $U\subset \p^3$ be an irreducible surface of degree 
$g\leq \sqrt[3]{6n}$.
Then there is a surface $T$ of degree $e$ 
that contains  ${\mathcal P}$,
does not contain $U$ and $ge^2\leq 6n$. 
\end{enumerate}
\end{lem}

Proof.  We argue as in Lemma \ref{lowder.surf.lem}.
For a surface of degree $d$ it is $1$ linear condition to
contain a point. Thus if $\binom{d+3}{3}>n$ then
such a surface $S$ exists, giving (1).

In order to prove (2) we need to find $e$ such that
$$
\tbinom{e+3}{3}-\tbinom{e-g+3}{3}> n.
$$
As before, the left had side is $\tfrac16 g(e+1)(e+5)>\tfrac16 g(e+1)^2$.
Thus  we can choose $e:=\rdown{\sqrt{6n/g}}$.
\qed

\begin{rem}  
Over infinite fields,  both lemmas can be extended to the case
when we want to avoid any finite collection of surfaces $U_i$
whose degrees are between $g$ and  $\sqrt{6m}$ in Lemma \ref{lowder.surf.lem}
(resp.\  between $g$ and  $\sqrt[3]{6n}$ in Lemma \ref{lowder.surf.pts.lem}).
We just need to take a general linear combination of the equations
obtained for the individual $U_i$.
\end{rem}

The conclusions of the second part of 
Lemmas \ref{lowder.surf.lem}--\ref{lowder.surf.pts.lem} get
weaker as $g$ gets smaller. I believe that Lemma \ref{lowder.surf.pts.lem}
 can not be  improved, but a
quite different method works for Lemma \ref{lowder.surf.lem}.

Let $S\subset \p^3$ be a surface of degree $d$. 
 In 1849 Salmon wrote down an equation of degree  $11d-24$
that cuts out on $S$ the locus of points where there is a
triple tangent line; see
\cite[pp.277--291]{salmon-2-old} for a detailed treatment
based on \cite{clebsch-61}. This locus clearly contains the union of
all lines contained in $S$. Cayley noted that
 every point has a triple tangent line iff $S$ is ruled.
The latter assertion is already in the fourth  edition of
Monge's book
\cite[\S XXI]{monge-old}, see especially p.225.
(I could not find the 1801 first edition
{\it Feuilles d'analyse appliqu\'ee \`a la g\'eom\'etrie;}
it is much shorter than the 1809 fourth  edition.)

\begin{thm}[Monge--Salmon--Cayley] \label{cayley-salmon.thm}
Let $S\subset \c\p^3$ be a surface  of degree $d$ 
without ruled irreducible components. 
Then there is a surface $T$ of degree $11d-24$ such that
$S$ and $T$ do not have common irreducible components and
every line on $S$ is contained in $S\cap T$. 
\end{thm}

For the reader's convenience, I give a---partly analytic---proof of
this in Section \ref{sketch.sec}. Strictly speaking, I only
show that $\deg T\leq 11d-18$. In the applications I use
only that $\deg T\leq 11d$, so this is not a problem.
For an algebraic approach  see \cite{MR1986126}, where the emphasis is
on understanding what happens in positive characteristic.

Once we have two surfaces $S,T$, we use the following
bound on the number of intersections.
This is the  observation that makes the estimates
in the Theorems     readily computable.

Let $C$ be a reduced curve. For a point $p\in C$, let
$r(p)$ denote the {\it multiplicity} of $C$ at $p$.
For line configurations, this equals the number of lines passing through $p$.
Since we use only the line configuration case, we do not discuss
the extra complications that appear in general.

\begin{prop} \label{ci.genus.prop}
Let $S, T\subset \p^3$ be two surfaces of degrees $a$
and $b$ that have no common irreducible components. 
Set $C=S\cap T$  (with reduced structure).  Then
\begin{enumerate}
\item $C$ has at most $ab$ irreducible components.
\item $\tsum_{p\in C} \ \bigl(r(p)-1\bigr)\leq \tfrac12 ab(a+b-2)$.
\item $\tsum_{p\in C} \ \bigl(r(p)-1\bigr)^{3/2}\leq \tfrac1{\sqrt{2}} ab(a+b-2)$.
\item $\tsum^{(sm)}_{p\in C}\  r(p)\bigl(r(p)-1\bigr)\leq  ab(a+b-2)$
where the sum is over those points where either $S$ or $T$ is smooth.
\end{enumerate}
\end{prop}

Outline of   proof.
We repeatedly use the {\it theorem of B\'ezout}
which says that if $H_1,\dots, H_n\subset \p^n$
are hypersurfaces of degrees $d_1,\dots, d_n$ then
their intersection $H_1\cap\cdots\cap H_n$ 
either contains an algebraic curve or it consist of at most $ d_1\cdots d_n$
points; cf. \cite[Sec.IV.2.1]{shaf}. 

Using this for $S,T$ and a general hyperplane we see that 
$C$ has degree $\leq ab$, thus
 $\leq ab$ irreducible components; if equality holds  then 
all irreducible components are lines, proving (\ref{ci.genus.prop}.1).

The proof of (\ref{ci.genus.prop}.2--4) has 2 main steps.

Note that $\tfrac12 ab(a+b-4)+1$ is the genus
of a smooth complete intersection curve of  two surfaces of degrees $a$
and $b$. This is a well known formula; see for example
\cite[Sec.VI.1.4]{shaf} (especially Exercise 9 on p.68 of volume 2)
or \cite[Exrc.I.7.2]{hartsh}.
The key claim is that even very singular 
complete intersection curves have arithmetic genus
$\leq \tfrac12 ab(a+b-4)+1$; see
Section \ref{sec.subcurves} for details.

Note that the  arithmetic genus frequently jumps up for
singular curves in families. 
(Historically,  schemes
and flatness were  introduced to understand similar phenomena.)
For instance, all rational curves of degree
$d$ in $\p^3$ form a single family. General members are smooth,
thus with genus $0$. At the other extreme we get plane rational curves
of degree $d$, these have arithmetic genus  $\binom{d-1}{2}$.

The second step is to use the arithmetic genus of a curve to
control its singularities and convert this information
into the estimates (\ref{ci.genus.prop}.2--4).
This is done in Section \ref{coh.sec}.

\begin{rem}\label{2.11.thm.rem}
 \cite{gu-ka} suggests (see especially Proposition 2.2 and the Appendix)
that, at least over $\r$,
for line configurations the optimal bound is of the form
$$
 \tsum^*_p \ \bigl(r(p)-1\bigr)^{2}\leq  (\mbox{constant})\cdot ab(a+b-2)
\log(a+b),
\eqno{(\ref{2.11.thm.rem}.1)}
$$
where summation is over the points satisfying   $1\leq r(p)\leq a+b$.
The appearance of $\log$ on the right hand side is  surprising
from the point of view of algebraic geometry.

I do not know if (\ref{2.11.thm.rem}.1)  holds over $\c$ or not, but
 over finite fields the exponent $3/2$ is optimal
as shown by  Example \ref{charp.many.intersections}.
The  exponent $3/2$ is also optimal for complete intersection curves
in general, even when the singularities locally look like
unions of lines. 

As an example, pick general homogeneous polynomials
$f,g,h$ of degree $n$. For general  $\alpha, \beta, \gamma\in \r$  set
$$
S:=\bigl(f^m+g^m+h^m=0\bigr)\qtq{and}
T:=\bigl(\alpha f^m+\beta g^m+\gamma h^m=0\bigr).
$$
Then $C:=S\cap T$ has $n^3$ singular points (where $f=g=h=0$)
and, at each of these points $C$ has $m^2$ smooth branches. Thus
$$
\tsum_{p\in C} \ \bigl(r(p)-1\bigr)^{3/2}=
n^3(m^2-1)^{3/2}
\leq
(nm)(nm)(2nm-2)
$$
indeed holds but the exponent $3/2$ can not be replaced with anything bigger.

We can even arrange all the singular points to be real.
\end{rem}

\begin{rem}  The maximum possible number of lines on a
degree $d$ non-ruled surface is not known. 
The  Fermat-type surface
$$
F_d:=\bigl(x_0^d+x_1^d+x_2^d+x_3^d=0\bigr)
$$
 contains $3d^2$ lines. There are a few examples with more lines,
for instance there are degree $20$ surfaces with  $4\cdot 20^2$ lines. 
See \cite{MR2341513, 2012arXiv1212.3511R, 2013arXiv1303.1304R} 
and the references there for further examples.
Over finite fields one can have many more lines, see
Example \ref{d4.lines.charp.exmp}.
\end{rem}

\section{Counting intersections}\label{sec.intersections}

\begin{say}\label{pf.1.setup.say}
 In order to prove Theorem \ref{main.thm},
let  $S$ be a  surface of smallest possible degree
that contains  our  $m$  lines $[{\mathcal L}]$.
By Lemma \ref{lowder.surf.lem} we know that $d:=\deg S\leq \sqrt{6m}$. 

Fix an ordering of the  irreducible
components $S_i\subset S$ and let  ${\mathcal L}_i\subset {\mathcal L}$
denote those lines that are contained in $S_i$ but are
not contained in $S_1\cup\cdots\cup S_{i-1}$. We can write
$I({\mathcal L}) $ in the form 
$$
I({\mathcal L})=\tsum_{i= 1}^r I({\mathcal L}_i)
+ \sum_{i= 2}^r \#\bigl(|{\mathcal L}_i| \cap 
|{\mathcal L}_1\cup\cdots \cup {\mathcal L}_{i-1}|\bigr).
\eqno{(\ref{pf.1.setup.say}.1)}
$$
The second sum counts  intersections of lines that lie on different
irreducible components; these are easy to bound, see (\ref{ext.say}).

For the first sum we need to work with one irreducible component
at a time. We treat 3 cases separately: planes and quadrics
in (\ref{pl+quad.say}), ruled surfaces of degree $\geq 3$ in (\ref{ruled.say})
and non-ruled surfaces in (\ref{better.cs.cor}).
\end{say}

\begin{say}[External intersections] \label{ext.say}
Let $\ell\in {\mathcal L}$
be a line and $S'(\ell)\subset S$ the union of those irreducible
components of $S$ that do not contain $\ell$. Any intersection
point of  $\ell$ with a line that is in $S'(\ell)$ is also
contained in $\ell\cap S'(\ell)$. By B\'ezout, this is  a set
of at most $\deg S'(\ell)\leq \deg S=d$ elements.
This gives at most $md$ such intersections.
We can do  better if we order the $S_i$ such that
$\#{\mathcal L}_i/\deg S_i$ is a non-increasing function of $i$.
(Note that  the set ${\mathcal L}_i$ depends on the ordering
of the surfaces, but we can choose at each step the surface that
maximizes the quotient.)
With such a choice there are at most $\frac12 md$ external intersections.
\end{say}



We are thus left to  work with  the surfaces $S_i$ separately
and  estimate the number of  intersections of the lines in  ${\mathcal L}_i$.

\begin{say}[Planes and quadrics] \label{pl+quad.say}
Let  $\{P_i: i\in I\}$  be the planes  and  $\{Q_j: j\in J\}$ the  quadrics
in $S$.
Set $m_i=\# {\mathcal L}_i$ for $i\in I$ 
and  $n_j=\# {\mathcal L}_j$ for $j\in J$.

A line on $P_i$ intersects all the other lines on
$P_i$ thus  $I({\mathcal L}_i)\leq \frac12 m_i^2$.
By assumption  $m_i\leq c\sqrt{m}$ 
thus  $I({\mathcal L}_i)\leq \frac12 m_ic\sqrt{m}$. 

On a singular quadric, any two lines meet at the vertex, thus
$I({\mathcal L}_j)= n_j-1$. On a smooth quadric, there are
2 families of lines. Correspondingly write $n_j=n'_j+n''_j$.
Then  
$$
I({\mathcal L}_i)\leq n'_jn''_j\leq \tfrac14 (n'_j+n''_j)^2=\tfrac14 n_j^2
\leq \tfrac14 n_j 2c\sqrt{m}.
$$
Thus the number of internal intersections on planes and quadrics
is at most
$$
\tsum_{i\in I}\tfrac12 m_ic\sqrt{m}+ \sum_{j\in J}\tfrac12 n_jc\sqrt{m}=
\tfrac12 
c\sqrt{m}\Bigl(\tsum_{i\in I}m_i+ \sum_{j\in J}n_j\Bigr)\leq \tfrac12 cm^{3/2}.
$$
\end{say}

\begin{say}[Ruled surfaces] \label{ruled.say}
Basic results on ruled surfaces are discussed in Section \ref{ruled.sec}.

Let $\{R_i:i\in I\}$
be the  irreducible ruled surface in $S$ and set $d_i=\deg R_i$.

By (\ref{ruled.surf.long.prop}.4) there are at most 2  lines, 
called {\it special lines,}
that intersect infinitely many  other lines. For each irreducible ruled surface
these contribute at most $2m_i$. By (\ref{ruled.surf.long.prop}.5),
every non-special line intersects at most $d_i:=\deg R_i$ 
other non-special  lines,
 hence these contribute at most $\frac12 m_id_i$.

Thus all together
we get at most $\frac12 md+2m$ intersections.
\end{say}

Finally we give two different bounds
for the non-ruled irreducible components.
First, combining Theorem \ref{cayley-salmon.thm} with
Proposition \ref{ci.genus.prop}.2 we get the following.

\begin{cor} \label{cayley-salmon.cor}
Let $S\subset \c\p^3$ be a surface of degree $d$ 
without ruled irreducible components and
 ${\mathcal L}$ the set of lines on $S$. Then
\begin{enumerate}
\item ${\mathcal L}$ contains at most $d(11d-24)$ lines and
\item  
$I({\mathcal L})\leq \tfrac12 d(11d-24)(12d-28)+ d(11d-24)\leq 66d^3$. \qed
\end{enumerate}
\end{cor}

The above bound does not involve $m$, so it is  best when
the degree of the surface $S$ is  small compared to
the number of lines. When $\deg S$ is close to the 
bound $\sqrt{6m}$ given in Lemma \ref{lowder.surf.lem}.1,
we get a better estimate using Lemma \ref{lowder.surf.lem}.2.

\begin{prop} \label{better.cs.prop}
 Let ${\mathcal L}$ be a set of $m$ distinct lines in $\p^3$
and $S\subset \p^3$  a minimal degree surface containing $[{\mathcal L}]$.
Assume that $S$ is  irreducible and has degree $d$.
 Then  
$I({\mathcal L})\leq 3m\bigl(d+\tfrac{6m}{d}\bigr)$.
\end{prop}

Proof. By Lemma \ref{lowder.surf.lem}, there is another surface $T$ of degree 
$\leq \tfrac{6m}{d}$ that contains $[{\mathcal L}]$. Applying
Proposition \ref{ci.genus.prop}.2 to $S,T$ we get our bound.\qed
\medskip

\begin{cor} \label{better.cs.cor}
 Let ${\mathcal L}$ be a set of $m$ distinct lines in $\p^3$
and $S\subset \p^3$  a minimal degree surface containing $[{\mathcal L}]$.
Assume that $S$ is  irreducible and   non-ruled.
 Then   
$I({\mathcal L})\leq 26.6\cdot m^{3/2}$.
\end{cor}

Proof. Set $d:=\deg S$ and
write it as  $d=\alpha \sqrt{m}$.  Note that
$\alpha\leq \sqrt{6}$  by Lemma \ref{lowder.surf.lem}.
Both Corollary  \ref{cayley-salmon.cor} and Proposition \ref{better.cs.prop}
give bounds, thus 
$$
I({\mathcal L})\leq 
\min\bigl\{ 66\alpha^3, 3\bigl(\alpha+\tfrac{6}{\alpha}\bigr)\bigr\}
\cdot m^{3/2}.
$$
The minimum reaches its maximum when the two quantities are equal.
This happens at
$\alpha_0 = \sqrt{6/11}\approx 0.738$
and  $66\alpha_0^3< 26.6$. \qed

\begin{say}[Adding up] Starting with $m$ distinct lines
${\mathcal L}$, let $S$ be the smallest degree surface
that contains $[{\mathcal L}]$. Note that each irreducible component 
$S_i\subset S$
has minimal degree among those surfaces that contain every line of 
${\mathcal L}_i$ (as in Paragraph \ref{ext.say}). 
We have 4 sources of intersection points.

External intersections  (\ref{ext.say}) contribute $\leq \frac12 md$,
planes and quadrics (\ref{pl+quad.say}) contribute 
$\leq  \frac12 cm^{3/2}$ and the other
ruled surfaces (\ref{ruled.say}) contribute $\leq \frac12 md+2m$.

Let $\{S_i:i\in I\}$ be the  non-ruled  irreducible components
and $m_i$ denote the number of lines in ${\mathcal L}_i$.
By Corollary \ref{better.cs.cor} these lines have at most
$26.6 m_i^{3/2}$ intersections with each other. Thus the
non-ruled  irreducible components contribute at most
$$
\tsum_{i\in I} 26.6 m_i^{3/2}\leq 26.6 m^{3/2}.
$$
So the total number of intersection points is at most
$$
  \tfrac12  md+ \tfrac12 cm^{3/2}+\tfrac12 md+2m +26.6 m^{3/2}
$$
Since $d\leq \sqrt{6m}-2$ by Lemma \ref{lowder.surf.lem}, this is at most
$$
\bigl(\sqrt{6}+26.6+\tfrac12 c\bigr)m^{3/2}
< \bigl(29.1+\tfrac{c}2 \bigr)m^{3/2}.
$$
This completes the proof of Theorem \ref{main.thm}.   \qed
\end{say}

\section{Counting incidences over $\c$}\label{sec.incidences.C}

\begin{say}\label{i.i.defn.say}
Let ${\mathcal L}$ be a set of $m$ distinct lines 
and ${\mathcal P}$  a set of $n$ distinct points
in $\p^3$.
Instead of $I({\mathcal L},{\mathcal P}) $
it is more convenient to work with the smaller quantity
$$
I^{\circ}({\mathcal L},{\mathcal P}):=\tsum_{p\in{|\mathcal L}|\cap{\mathcal P} } 
\bigl(r(p)-1\bigr)
\eqno{(\ref{i.i.defn.say}.1)}
$$
which is  better suited to
induction thanks to the  subadditivity property:
$$
I^{\circ}({\mathcal L}\cup\{\ell\},{\mathcal P})\leq
I^{\circ}({\mathcal L},{\mathcal P})
+ \#\bigl([{\mathcal L}]\cap\ell\bigr)
\qtq{provided} \ell\notin {\mathcal L}.
\eqno{(\ref{i.i.defn.say}.2)} 
$$
The two variants are related by the formula 
$I({\mathcal L}, {\mathcal P})=I^{\circ}({\mathcal L},{\mathcal P})
+\# \bigl([{\mathcal L}]\cap{\mathcal P}\bigr)$.
\end{say}

As a preliminary step toward proving Theorem \ref{main.thm.2}
 we reduce to the case when
every line meets ${\mathcal P} $ in many points.

\begin{say}[Lines with few points]\label{few.pts.l.say}
Assume that under the assumptions of  Theorems \ref{main.thm.2} 
or \ref{main.thm.2.fq} 
 we want to prove a bound  of the form 
$$
I({\mathcal L},{\mathcal P})\leq   
 m A(n)+(c^2m) B(n)+C(n)
\eqno{(\ref{few.pts.l.say}.1)}
$$
for some functions $A(n), B(n), C(n)$.
Let $\ell\in {\mathcal L}$ be a line 
 that  meets ${\mathcal P} $ in 
 $\leq A(n)$ points. Remove  $\ell$ from ${\mathcal L}$.
 Note that we may need to
increase  $c$  to $c(\frac{m}{m-1})^{1/2}$.  Thus the left hand side
of (\ref{few.pts.l.say}.1)
decreases by  $\leq A(n)$ and the right hand side by
$$
m\bigl(A(n)+c^2B(n)\bigr) -
(m-1)\bigl(A(n)+c^2\tfrac{m}{m-1}B(n)\bigr) =A(n).
$$
Hence it is sufficient to prove (\ref{few.pts.l.say}.1) for
line/point configurations where every line meets ${\mathcal P} $ in 
 $> A(n)$ points.

This step makes the proof less direct. In
Section \ref{sec.intersections} we just wrote down the estimates and
got a final result. Here  we need to know in advance the final result
we aim at and use the corresponding value of $A(n)$.
\end{say}

\begin{say}[Decomposing $S$ and ${\mathcal P}$] \label{start.of.pf.thm2}
Let  $S$ be a  surface of smallest possible degree
that contains  our set of $n$ distinct points ${\mathcal P}$.
By Lemma \ref{lowder.surf.pts.lem} we know that $d:=\deg S\leq \sqrt[3]{6n}$.

We would like to ensure that
$S$ contains all the lines in  ${\mathcal L}$.
If a line $\ell$ is not contained in $S$ then, by B\'ezout, 
it meets $S$ in at most  $d\leq \sqrt[3]{6n}$ points.
Thus if $\ell$ passes through more than $\sqrt[3]{6n}$ points
of ${\mathcal P}$ then  $\ell\subset S$. 
This suggests that we use (\ref{few.pts.l.say}) with
 $A(n)=1.82n^{1/3}>\sqrt[3]{6}n^{1/3}$.  Thus we may assume that 
 each line  in 
${\mathcal L}$ contains $\geq  1.82n^{1/3}$ points
of  $ {\mathcal P} $ hence
$[{\mathcal L}]$ is contained in $S$.

We will also need to
divide the points among the  irreducible components of $S$.
Let  $S_i\subset S$ be an irreducible component of degree $d_i$.
Let ${\mathcal P}^*_i\subset{\mathcal P}$ denote the subset of points
 that are  on $S_i$ but not on any other
irreducible component of $S$. There is at most 1 component, call it $S_0$,
for which $|{\mathcal P}^*_0|>\frac12 n$. 
Let  ${\mathcal P}_0\subset{\mathcal P}$ denote the subset of points
 that are  on $S_0$; for $i\neq 0$ set ${\mathcal P}_i={\mathcal P}^*_i$.
The ${\mathcal P}_i$ are disjoint subsets of ${\mathcal P}$,
thus $\sum n_i\leq n$ where  $n_i:=|{\mathcal P}_i|$.
 Since $S$ has minimal degree, we know that
each $S_i$ has minimal degree among those surfaces that contain
${\mathcal P}_i$, hence $d_i\leq \sqrt[3]{6n_i}$.

Next we use $S$ to estimate $ I^{\circ}({\mathcal L},{\mathcal P})$.
As before
we try to find another surface $T$ that contains ${\mathcal L} $
but does not contain $S$ or at least
some of the irreducible components of $S$.
\end{say}

\begin{say}[Contributions from singular points of $S$]
\label{pf.of.main.thm.2.3}
We start with  lines contained in $\sing S$ and their 
intersection points. If $S$ is defined by an equation 
$\bigl(f(x_0,\dots,x_3)=0\bigr)$ then we can take $T$
to be defined by a general linear combination
$$
\tsum_i a_i\frac{\partial f}{\partial x_i}=0.
$$
Thus $\deg T=d-1$ and, using (\ref{ci.genus.prop}.2),
we get  a contribution to  $ I^{\circ}({\mathcal L},{\mathcal P})$
that is  $\leq \tfrac12 d(d-1)(2d-3)\leq d^3\leq 6n$.
 This is the contribution 
 from  lines that are contained in $\sing S$. 

We can do  better using (\ref{ci.genus.prop}.3)
which says that
$$
\tsum\ \bigl(r(p)-1\bigr)^{3/2}\leq  6\sqrt{2}\cdot n.
\eqno{(\ref{pf.of.main.thm.2.3}.1)}
$$
Since we have at most $n$ summands on the left, 
the convexity of  $x^{3/2}$ implies that 
$$
\tsum\ \bigl(r(p)-1\bigr)\leq (6\sqrt{2})^{2/3}\cdot n< 4.2n.
\eqno{(\ref{pf.of.main.thm.2.3}.2)}
$$

Now we add to this lines $\ell_i$ not contained in $\sing S$ 
one at a time. Each  line intersects $\sing S$ in  at most $d-1$ points.
Repeatedly using (\ref{i.i.defn.say}.2)
we get a contribution of $\leq m(d-1)$.

These two account for all the contributions in
(\ref{i.i.defn.say}.1) coming from
those points of ${\mathcal P}$ that are singular on $S$.
\end{say}

\begin{say}[Contributions from smooth points of $S$; ruled case]
\label{pf.of.main.thm.2.4}
A smooth point is contained in a unique  irreducible component of $S$,
thus we can treat the irreducible components  separately.
We start with the ruled components.

\ref{pf.of.main.thm.2.4}.1 {\it Planes.}
By assumption, each plane contributes $\leq \tfrac12 c^2m$. Since
there are $\leq d$ planes,
 all together they contribute $\leq \tfrac12 c^2md$.

\ref{pf.of.main.thm.2.4}.2  {\it Other ruled surfaces.}
On a smooth quadric, there are 2 lines through each point.
On other ruled surfaces there is usually only 1 line
through a smooth point, except on the special lines (\ref{ruled.specv.say})
when there can be 2 by (\ref{ruled.surf.long.prop}).
Thus  we get  a  total contribution
$\leq n$.
\end{say}

As in Section \ref{sec.intersections}, we again  use 2 methods to 
control non-ruled irreducible components.

\begin{say}[Contributions from smooth points of $S$; non-ruled case I]
\label{pf.of.main.thm.2.5}
Let  $S_i\subset S$ be a non-ruled irreducible component of degree $d_i$.
As we noted in (\ref{start.of.pf.thm2}),  $ d_i\leq \sqrt[3]{6n_i}$.

By Theorem \ref{cayley-salmon.thm} there is another surface $T_i$ of degree
$\leq 11d_i$ that contains every line lying on $S_i$.
Using  (\ref{ci.genus.prop}.4) we get that
$$
\tsum^{(sm)}_i\ \bigl(r(p)-1\bigr)^2\leq 11\cdot 12\cdot d_i^3,
\eqno{(\ref{pf.of.main.thm.2.5}.1)}
$$
where summation is over all 
smooth points of $S_i$ that are in
${\mathcal P}_i$.
Since $d_i^3\leq  6 n_i$
and $\sum n_i\leq n$, adding these up gives that 
$$
\tsum^{(sm)}_S\ \bigl(r(p)-1\bigr)^2
\leq 11\cdot 12\cdot 6\cdot n.
$$
where summation is over all smooth points of the 
non-ruled irreducible components of $S$ that are in
${\mathcal P}$.
Since we have at most $n$ summands on the left, by Cauchy--Schwartz
$$
\tsum^{(sm)}_S\ \bigl(r(p)-1\bigr)
\leq \sqrt{11\cdot 12\cdot 6}\cdot n< 28.2\cdot  n.
\eqno{(\ref{pf.of.main.thm.2.5}.2)}
$$
\end{say}

\begin{say}[First estimate]\label{pf.of.main.thm.2.66}
Adding these together we get that
$$
\begin{array}{rcl}
I({\mathcal L}, {\mathcal P})&\leq &
n+I^{\circ}({\mathcal L},{\mathcal P})\\
&\leq & n+4.2n+m(d-1)+\tfrac12 c^2md+n+28.2  n\\
&\leq & \bigl(1+\tfrac12 c^2\bigr)md+34.6  n\\
&\leq & \sqrt[3]{6}\bigl(1+\tfrac12 c^2\bigr)
mn^{1/3}+34.6n. 
\end{array}
\eqno{(\ref{pf.of.main.thm.2.66}.1)}
$$
This is different from the  bound claimed in Theorem \ref{main.thm.2}.
The coefficient of $mn^{1/3} $ is  smaller but the
the coefficient of $n $ is  bigger. 
For some applications this may be a better variant
but (\ref{pf.of.main.thm.2.66}.1) gives a 
worse constant for Corollary \ref{main.thm.bourg.2}.
\end{say}

We need to look at the non-ruled components again.

\begin{say}[Contributions from smooth points of $S$; non-ruled case II]
\label{pf.of.main.thm.2.6}
Here we are aiming to get an estimate as in (\ref{few.pts.l.say}.1) with
$A(n)=3.66n^{1/3}$ which is chosen to be an upper bound for
 $ \sqrt{6\sqrt[3]{11}}n^{1/3}$. 
Thus we may assume that each line contains
$\geq 3.66n^{1/3}$ points of ${\mathcal P}$.

Write  $d_i=\alpha_i\cdot n_i^{1/3}$ and note that $\alpha_i\leq \sqrt[3]{6}$. 
We improve the previous estimate if
$\alpha_i\geq 1/\sqrt[3]{11}$. 
By Lemma \ref{lowder.surf.pts.lem} there is a surface $T_i$ 
of degree 
$\leq \sqrt{6n_i/d_i}= \sqrt{6/\alpha_i}\cdot n_i^{1/3}$
that contains ${\mathcal P}_i$ but not $S_i$.

If $i=0$ then every line in ${\mathcal L}$
that is contained in $S_0$ meets ${\mathcal P}$, and hence
also  ${\mathcal P}_0={\mathcal P}\cap S_0$, in at least
$3.66n^{1/3}$ points. Thus these lines  are also contained in $T_0$.

If $i>0$ then let  $T^{(i)}$ be the surface obtained from $S$ by
replacing $S_i$ with $T_i$. Note that  $T^{(i)}$  contains ${\mathcal P}$ and
its degree is 
$$
\leq \sqrt{6/\alpha_i}n_i^{1/3}+d-d_i\leq
\bigl(\sqrt{6/\alpha_i}-\alpha_i\bigr)n_i^{1/3}+ \sqrt[3]{6}n^{1/3}.
$$
 Since $n_i<\frac12 n$, this is less than
$3.66n^{1/3}$. Thus  $T^{(i)}$ contains ${\mathcal L}$ and hence
 $T_i$ contains every line in ${\mathcal L}$
that is not contained in any other $S_j$. 

Since $\alpha_i\leq \sqrt[3]{6}$, this gives a bound 
$$
\begin{array}{rcl}
\tsum^{(sm)}_i\ \bigl(r(p)-1\bigr)^2&\leq& 
\alpha_in_i^{1/3}\sqrt{6/\alpha_i}n_i^{1/3}
\bigl(\alpha_in_i^{1/3}+\sqrt{6/\alpha_i}n_i^{1/3}\bigr)\\
& =& \bigl(6+\sqrt{6}\alpha_i^{3/2}\bigr)n_i\leq 12 n_i.
\end{array}
\eqno{(\ref{pf.of.main.thm.2.6}.1)}
$$
If $\alpha_i\leq 1/\sqrt[3]{11}$ then
 $d_i\leq \bigl(1/\sqrt[3]{11}\bigr)n_i^{1/3}$ 
and so (\ref{pf.of.main.thm.2.5}.1) and (\ref{pf.of.main.thm.2.6}.1)
together show that 
$$
\tsum^{(sm)}_i\ \bigl(r(p)-1\bigr)^2\leq  12n_i
\eqno{(\ref{pf.of.main.thm.2.6}.2)}
$$
holds for every non-ruled surface.
Adding up all cases gives that 
$$
\tsum^{(sm)}_S\ \bigl(r(p)-1\bigr)^2\leq 12n.
\eqno{(\ref{pf.of.main.thm.2.6}.3)}
$$
As before,  by Cauchy--Schwartz this implies that 
$\tsum^{(sm)}_S\  \bigl(r(p)-1\bigr)\leq \sqrt{12}  n$.
\end{say}

\begin{say}[Final estimate II]\label{pf.of.main.thm.2.8}
Adding these together we get that
$$
\begin{array}{rcl}
I({\mathcal L}, {\mathcal P})&\leq &
n+I^{\circ}({\mathcal L},{\mathcal P})\\
&\leq & n+2\cdot 3^{2/3}n+m(d-1)+\tfrac12 c^2md+n+12^{1/2}  n\\
&\leq & \bigl(1+\tfrac12 c^2\bigr)md+9.9  n\\
&\leq & \sqrt[3]{6}\bigl(1+\tfrac12 c^2\bigr)
mn^{1/3}+9.9n. 
\end{array}
\eqno{(\ref{pf.of.main.thm.2.8}.1)}
$$
There is one place where it is easy to improve the estimate.
Assume that there are $xn$ points  used in 
(\ref{pf.of.main.thm.2.3}.1), $yn$ points used in 
(\ref{pf.of.main.thm.2.4}.2) and $zn$ points used in 
(\ref{pf.of.main.thm.2.6}.3). Then
$x+y+z\leq 1$ and the total contribution coming from
these points is at most
$$
\bigl(2\cdot 3^{2/3}x^{1/3}+y+2\cdot 3^{1/2}z^{1/2}\bigr)\cdot n
\qtq{where} x+y+z\leq 1.
$$
A straightforward computation using Lagrange multipliers
shows that this is always $\leq 5.76n$.
Thus we get that
$$
I({\mathcal L}, {\mathcal P})\leq 
\sqrt[3]{6}\bigl(1+\tfrac12 c^2\bigr)mn^{1/3}+6.76n.
\eqno{(\ref{pf.of.main.thm.2.8}.2)}
$$
Note however, that we assumed that each line contains
$\geq 3.66n^{1/3}$ points of ${\mathcal P}$.
By the first reduction step
(\ref{few.pts.l.say}) this  requires  us to have  $A(n)\geq  3.66n^{1/3} $
in (\ref{few.pts.l.say}.1), thus we can not use the smaller coefficient 
$\sqrt[3]{6}\leq 1.82$ in general. 
Thus we can only conclude that
$$
I({\mathcal L}, {\mathcal P})\leq 
\max\{3.66, 1.82+0.91c^2\}mn^{1/3}+6.76n.
\eqno{(\ref{pf.of.main.thm.2.8}.3)}
$$
This is stronger than Theorem \ref{main.thm.2}. \qed
\end{say}

\begin{say}[Proof of Corollary \ref{main.thm.bourg.2}]
\label{pf.main.thm.bourg.2}
We start with the characteristic 0 case.
Note that $50>3.66^3$, thus 
if $n<\tfrac1{50}m^{3/2}$ then $\sqrt{m}>3.66 n^{1/3}$.
This means that we do not need to go through the first reduction step
(\ref{few.pts.l.say}), hence the stronger conclusion
(\ref{pf.of.main.thm.2.8}.2) applies. 

 Choose $x$ such that   $n=\tfrac1{x^3}m^{3/2}$.
Since we assume that  $c=1$, (\ref{pf.of.main.thm.2.8}.2) becomes
$$
I({\mathcal L}, {\mathcal P})\leq
\bigl({2.73} \tfrac1{x}+6.76\tfrac1{x^3}\bigr)m^{3/2}.
$$
We compute that if $x^3\geq 50$ then 
${2.73}\tfrac1{x}+6.76\tfrac1{x^3} <1$ hence
$I({\mathcal L}, {\mathcal P})<m^{3/2}$.
On the other hand, by assumption each line contributes at
least $m^{1/2}$, hence $I({\mathcal L}, {\mathcal P})\geq m^{3/2}$.
This is a contradiction if $n< \tfrac1{50}\cdot m^{3/2} $.

The positive  characteristic case follows from Theorem \ref{main.thm.2.fq}
similarly. For large $m$ the proof gives a coefficient  $\geq \tfrac1{13}$;
the smaller value $\tfrac1{20}$ and the 
$m\geq 10^4$ assumption are there to account for the  contribution of
 the two lower
degree terms in  (\ref{main.thm.2.fq}.1).
\qed
\end{say}


\section{Counting incidences over $\f_q$}\label{sec.incidences.Fq}

In this section we work with arbitrary fields, but the main interest
is understanding what happens over finite fields.

While Salmon's argument applies over any field,
Monge's proof only works in characteristic 0. 
As a replacement, (\ref{lines.on.surf.prop})
shows that a surface of degree $d$ 
without ruled irreducible components contains at most $d^4$ lines.
The following example shows that this is quite sharp.

\begin{exmp} \label{d4.lines.charp.exmp}
Let $q$ be a $p$-power and consider the surface
$$
S_{q+1}:=\bigl(x^{q+1}+y^{q+1}+z^{q+1}+w^{q+1}=0\bigr)\subset\p^3
$$
over the field $\f_q$. Linear spaces on such  {\it Hermitian
hypersurfaces} have been studied in detail
\cite{MR0213949, MR0200782}. 
These examples have also long been recognized as extremal
for the  Gauss map. The failure of
Monge's theorem has been noted in \cite{MR1143550} for surfaces
and in \cite{MR0080354} for curves. See 
\cite{MR846021, MR2431661} for  surveys of the Gauss map.
Other  extremal properties are discussed in \cite{2013arXiv1304.0302H}.
Kleiman observed that the affine Heisenberg surface in \cite[\S
8]{MR2031165} is, after taking its closure in $\p^3$, a Hermitian surface,
so isomorphic, under an  $\f_q$-linear transformation, to the surface above.

The configuration of lines on  $S_{q+1} $ is quite 
interesting.
\begin{enumerate}
\item  $S_{q+1} $ contains  $(q+1)(q^3+1) $ lines, all defined over $ \f_{q^2} $.
\item  $S_{q+1} $ contains  $(q^2+1)(q^3+1) $   points in  $ \f_{q^2} $.
\item $\PSU_4(q)$ acts transitively on the lines and on the  $ \f_{q^2} $-points.
\item There are $q+1$ lines through every  $ \f_{q^2} $-point.
\end{enumerate}
All of these are easy to do by hand as in  see \cite{MR0213949, MR0200782}
or can be obtained  from the general description of finite unitary groups;
see for instance \cite{MR0407163}.

More generally consider any equation of the form
$$
\tsum_{0\leq i,j\leq n} c_{ij} x_i^qx_j=0.
\eqno{(\ref{d4.lines.charp.exmp}.5)}
$$
If we substitute  $x_i=a_it+b_is$ then we get
$$
\tsum_{ij} c_{ij} (a_it+b_is)^q(a_jt+b_js)=
\tsum_{ij} c_{ij} (a_i^qt^q+b_i^qs^q)(a_jt+b_js)=0,
$$
which involves only the monomials  $t^{q+1}, t^qs, ts^q, s^{q+1}$.
Thus, arguing as in (\ref{lines.on.surf.prop}), we expect 
many more lines than usual.
It was proved by \cite{MR1550551} that if the hypersurface given by
(\ref{d4.lines.charp.exmp}.5) is smooth then it is isomorphic to
the Hermitian example, though the coordinate change is usually defined
 only over a field extension of $\f_q$.

\end{exmp}


The  arguments  in Section \ref{sec.incidences.C} are independent of the 
characteristic, save  
(\ref{pf.of.main.thm.2.5}) where 
we started considering non-ruled irreducible components.
We show below how to modify the estimates in
(\ref{pf.of.main.thm.2.5}--\ref{pf.of.main.thm.2.8}) to work over any field.

\begin{say}[Contributions from smooth points of $S$; non-ruled case I]
\label{pf.of.main.thm.2.5.fq}
Let  $S_i\subset S$ be a non-ruled irreducible component of degree $d_i$
and ${\mathcal P}_i\subset {\mathcal P}$ as in (\ref{start.of.pf.thm2}).

By Proposition \ref{lines.on.surf.prop} and 
Lemma \ref{lowder.surf.lem} there is another surface $T_i$ of degree
$\sqrt{6}d_i^3$ that contains every line lying on $S_i$.
Using  (\ref{ci.genus.prop}.4) we get that
$$
\tsum^{(sm)}_i\ \bigl(r(p)-1\bigr)^2\leq  6d_i^7+\sqrt{6}d_i^5,
\eqno{(\ref{pf.of.main.thm.2.5.fq}.1)}
$$
where summation is over all smooth points of $S_i$ 
that are in  ${\mathcal P}_i$.

Assume for now that   $ d_i\leq n_i^{1/5}$. Then
$6d_i^7+\sqrt{6}d_i^5\leq 6n_i^{7/5}+ \sqrt{6}n_i$.
Since we have at most $n_i$ summands on the left, by Cauchy--Schwartz
$$
\tsum^{(sm)}_S\ \bigl(r(p)-1\bigr)
\leq \sqrt{6}n_i^{6/5}+\tfrac12 n_i.
\eqno{(\ref{pf.of.main.thm.2.5.fq}.2)}
$$
\end{say}

\begin{say}[Contributions from smooth points of $S$; non-ruled case II]
\label{pf.of.main.thm.2.6.fq}
Here we deal with the other possibility   $ d_i\geq n_i^{1/5}$
 using (\ref{few.pts.l.say}) with
$A(n)=\sqrt{6}n^{2/5}$.

By Lemma \ref{lowder.surf.pts.lem} there is a surface $T_i$ 
of degree 
$\leq \sqrt{6n_i/d_i}\leq \sqrt{6}n_i^{2/5}$
that contains ${\mathcal P}_i$. If $i=0$ then 
$T_0$ contains every line of ${\mathcal L}$ that lies only on $S_0$.
 
If $i>0$ then $n_i\leq \frac12 n$ and, as in (\ref{pf.of.main.thm.2.6}),
 we get  a surface  $T^{(i)}$  of degree 
$$
\leq \sqrt{6}n_i^{2/5}+d-d_i< \sqrt{6}n_i^{2/5}+\sqrt[3]{6}n^{1/3}
< \sqrt{6}n^{2/5}
\eqno{(\ref{pf.of.main.thm.2.6.fq}.2)}
$$
that contains ${\mathcal P}$. 
 Therefore again
$T_i$ contains every line of ${\mathcal L}$ that lies only on $S_i$. 
The rest of  (\ref{pf.of.main.thm.2.6}) works as before and we get that
$$
\tsum^{(sm)}_i\  \bigl(r(p)-1\bigr)\leq \sqrt{12}  n_i.
\eqno{(\ref{pf.of.main.thm.2.6.fq}.3)}
$$
For $n_i\geq 2$ the right hand side of (\ref{pf.of.main.thm.2.5.fq}.2)
is bigger that $\sqrt{12}  n_i $, thus we can always use 
(\ref{pf.of.main.thm.2.5.fq}.2). \end{say}

\begin{say}[Final estimate]\label{pf.of.main.thm.2.fq}
Adding these together we get that
$$
\begin{array}{rcl}
I({\mathcal L}, {\mathcal P})&\leq &
n+I^{\circ}({\mathcal L},{\mathcal P})\\
&\leq & n+m(d-1)+\tfrac12 c^2md+ \\
&& 12^{2/3}n+n+\sqrt{6}n^{6/5}+\tfrac12 n.
\end{array}
\eqno{(\ref{pf.of.main.thm.2.fq}.1)}
$$
Note, however, that we have used (\ref{few.pts.l.say}) with
$A(n)=\sqrt{6}n^{2/5}$, thus the leading term $\sqrt[3]{6}mn^{1/3} $
needs to be increased to $\sqrt{6}mn^{2/5} $,
resulting in the final estimate
$$
\begin{array}{rcl}
I({\mathcal L}, {\mathcal P})
\leq 2.45mn^{2/5}+2.45n^{6/5} +0.91c^2mn^{1/3}+6.74n. 
\end{array}
\eqno{(\ref{pf.of.main.thm.2.fq}.2)}
$$
This completes the proof of (\ref{main.thm.2.fq}). 
As in (\ref{pf.of.main.thm.2.8})
we could improve the coefficient of $n$  a little but I see no
immediate application for it. \qed
\end{say}

\begin{say}[Bourgain's conjecture over finite fields]
\label{39.say}
We prove in Corollary \ref{main.thm.char.p}
that  Theorem \ref{main.thm}
holds over a field of characteristic $p$ provided
$p>\sqrt{m}$.
This implies that Theorem \ref{main.thm} holds for
all line configurations in $\f_p\p^3$ where $p$ is a prime.
For $p<\sqrt{m}$ the methods seem to yield only a weaker variant
with exponent $7/4$.

Similarly, Theorem \ref{main.thm.2} 
 holds
in characteristic $p$ provided $p>\sqrt[3]{6n}$. 
If we work over $\f_q$ then  
$I({\mathcal L}, {\mathcal P})\leq q^3+q^2+q+1$,
hence the estimate (\ref{main.thm.2}.1) is obvious 
if $q+1\leq  \sqrt[3]{6.76n}$.
Thus Theorem \ref{main.thm.2} holds over $\f_p$. 
(Note that  \cite{2013arXiv1311.1479E} gives
counter examples  over $\f_{p^2}$, building on \cite{MR2031165}.) 
\end{say}

The key to these is that Monge's theorem holds in characteristic $p>0$ 
 if the degree is
less than the characteristic. 
 \cite[Thm.1]{MR1986126} proves this for smooth surfaces
but essentially the same argument works in general.


\begin{cor}  \label{main.thm.char.p}
Let ${\mathcal L}$ be a set of $m$ distinct lines in $\f_q\p^3$ where $q=p^a$.
Let $c$ be a constant such that no plane (resp.\ no quadric) contains
more than $c\sqrt{m}$ (resp.\ more than $2c\sqrt{m}$) of the lines.

Assume that   either  $m<\frac{11}{6}p^2$ or $q=p$. 
Then the number of  points where at least two of the lines in
  ${\mathcal L}$ meet is  
$\leq (29.1+\frac{c}2)\cdot m^{3/2}$.
\end{cor}

Proof.  In the proof of Theorem \ref{main.thm}
 we used  Theorem \ref{cayley-salmon.thm}
only during the proof of  Corollary
\ref{better.cs.cor} where we applied it to a surface of degree 
$\leq \alpha_0\sqrt{m}$ with
$\alpha_0=\sqrt{6/11}$.
If  $m<\frac{11}{6}p^2$ then $\alpha_0\sqrt{m}<p$ 
hence, as noted above, 
Theorem \ref{cayley-salmon.thm} still
applies. 

If $q=p$ then we are done if $m<\frac{11}{6}p^2=\frac{11}{6}q^2$. 
If $m\geq \frac{11}{6}q^2$ then
we are done trivially since  $\f_q\p^3$ has $q^3+q^2+q+1$ points, hence
there are at most $2 m^{3/2}$ possible intersection points. \qed

\begin{exmp} \label{charp.many.intersections}
Let $L_1, L_2\subset \p^3$ be a pair of skew lines. For every point
$p\in \p^3\setminus(L_1\cup L_2)$ there is a unique line
$\ell_p$ passing through $p$ that intersects both  $L_1, L_2$.

The picture becomes especially simple when we work over a field $K$
and $L_1, L_2$ is a conjugate pair defined over a quadratic extension $K'/K$.
Thus we get that $\p^3(K)$ is a disjoint union of lines
naturally parametrized by the $K'$-points of $L_1$.
If $P\subset \p^3$ is a $K$-plane then $P\cap (L_1\cup L_2)$
consists of 2 points; the line connecting them is the only
line in our family that is contained in $P$.

For $K=\f_q$ we get a family of $q^2+1$ disjoint lines
$\{\ell_i\}$ that cover $\f_q\p^3$. 

A different pair of  skew lines  $L'_1, L'_2$
gives a different covering family of
lines $\{\ell'_i\}$. If  $L_1, L_2,L'_1, L'_2$ do not lie on a quadric
surface then they have $\leq 2$ common transversals.
(These are sometimes $K$-lines, sometimes conjugate pairs.)

Thus if we have $r$ different pairs of  skew lines
in general position then their union gives a family of
$m$ lines where 
$$
r(q^2+1)\geq m\geq r(q^2+1)-2\tbinom{r}{2}.
$$
The number of points where $r$ lines meet is at least
$$
q^3+q^2+q+1 -2(q+1)\tbinom{r}{2}.
$$
 Thus, for $r\ll\sqrt{m}$ we have
$$
m \qtq{lines and} \asymp \frac{m^{3/2}}{r^{3/2}}
\qtq{$r$-fold intersections.}
$$
Furthermore, any plane contains at most $r$ of the lines.

Given any set of $rq^2$ lines in $\f_q\p^3$, in average $r$ of them
pass through a point and $r$ of them are contained in a plane.
 The interesting aspect of the example is that for both of these,
the expected value is the maximum.
\end{exmp}

All of these examples either cover a positive proportion
 of $\f_q\p^3$ or can be derived 
by a linear transformation from a
configuration defined over  a subfield of
$\f_q$. It would be interesting if these turned out to be
the only cases that behave differently from characteristic 0.

\section{Genus and singular points of curves}\label{coh.sec}

\begin{say}[Hilbert polynomials] \label{hilb.say}
See  \cite[Chap.11]{at-mc} or \cite[Sec.I.7]{hartsh} for
proofs of the following results.

Let $k$ be a field, $R:=k[x_0,\dots, x_n]$ and
$I\subset R$  a homogeneous ideal.
The quotient ring $R/I$ is graded, that is,
it is the direct sum of its
homogeneous pieces $(R/I)_d$. Hilbert proved that there 
is a polynomial  $H_{R/I}(t)$, called the
{\it  Hilbert polynomial} of $R/I$ such that
$$
\dim (R/I)_d=H_{R/I}(d) \qtq{for $d\gg 1$.}
\eqno{(\ref{hilb.say}.1)}
$$
If $X\subset \p^n$ is a closed algebraic subvariety and
$I(X)$ the ideal of homogeneous polynomials that vanish on $X$ then
 $H_{R/I(X)}(t)$ is also called the
{\it  Hilbert polynomial} of $X$ and denoted by $H_X(t)$.

The degree of $H_{R/I}(t)$ equals the dimension of the
corresponding variety $V(I)$ and the leading coefficient
of $H_X(t)$ equals $\deg X/(\dim X!)$. The constant coefficient is
the (holomorphic) {\it Euler characteristic} of $X$,
denoted by $\chi(X, \o_X)$.

If $\dim X=1$ then, for historical reasons,
one usually uses the {\it arithmetic genus} 
$p_a(X):=1-\chi(X, \o_X)$.
(If $X$ is a smooth curve over $\c$ (=Riemann surface), the
arithmetic genus equals the topological genus.)

Let  $g\in k[x_0,\dots, x_n]$ be  homogeneous of degree $a$
and set $H:=(g=0)$. 
It is easy to see that if $g$
is not a zero-divisor on $X$ then
$$
H_{X\cap H}(t)=H_X(t)-H_X(t-a).
\eqno{(\ref{hilb.say}.2)}
$$
Assume next that we have hypersurfaces $H_i:=(g_i=0)\subset \p^n$ 
of degree $a_i$
such that $B:=H_1\cap \cdots\cap H_{n-1}$  has dimension 1.
(Such a $B$ is called a {\it complete intersection curve.})
Starting with 
$$
H_{\p^n}(t)=\tbinom{t+n}{n},
\eqno{(\ref{hilb.say}.3)}
$$
and using (\ref{hilb.say}.2) one can  compute the Hilbert polynomial of 
$B$:
$$
H_B(t)= \tprod_i a_i\cdot t-
\tfrac12 \bigl(\tsum_i a_i-n-1\bigr)\cdot \tprod_i a_i;
\eqno{(\ref{hilb.say}.4)}
$$
see  \cite[Exrc.II.8.4]{hartsh}. Thus the  arithmetic genus of $B$ is
$$
p_a(B):=1+\tfrac12 \bigl(\tsum_i a_i-n-1\bigr)\cdot \tprod_i a_i.
\eqno{(\ref{hilb.say}.5)}
$$

The formulas (\ref{hilb.say}.4--5) compute the Hilbert polynomial
and the arithmetic genus scheme-theoretically, that is,
we work with the Hilbert polynomial of the quotient ring
$k[x_0,\dots, x_n]/(g_1,\dots, g_{n-1})$
and this ring may contain nilpotents. 

As a simple example,  
consider  $B=(xy-zt=0)\cap (x(x+y)-zt=0)$, the intersection of two 
hyperboloids. Then $x\in k[x,y,z,t]/( xy-zt, x(x+y)-zt)$
is non-zero yet $x^2\in ( xy-zt, x(x+y)-zt)$.
The geometric picture is that $B$ consists of 
 2 lines  $L_1\cup L_2=(x=z=0)\cup (x=t=0)$, but $B$ ``counts'' both
 with multiplicity 2. The ideal corresponding to  $L_1\cup L_2$ is
$I(L_1\cup L_2)=(x, zt)$.

Thus the ideals $ ( xy-zt, x(x+y)-zt)$ and $(x, zt)$
define the same algebraic set. Given an algebraic curve
$B\subset \p^n$, it is usually not hard to find some of the equations
satisfied by $B$ and to write down an ideal $J$ that defines
$B$ set-theoretically.  However, we can not compute the arithmetic 
genus of $B$ using $J$.

It is usually much harder to write down
the ideal $I(B)$ of  {\em all} equations satisfied by $B$.
\end{say}

We prove the following basic inequality in the next section.

\begin{prop} \label{ag.prop.2}
For $i=1,\dots, n-1$ let $H_i\subset \p^n$ be 
a hypersurface of degree $a_i$ such that the  intersection
$B:=H_1\cap \cdots\cap H_{n-1}$ is 1-dimensional.
Let $C\subset B$ be a reduced subcurve. Then
 $$
p_a(C)\leq p_a(B)= 1+\tfrac12 \bigl(\tsum_i a_i-n-1\bigr)\cdot \tprod_i a_i .
$$
\end{prop}

\begin{say}[Arithmetic genus of a union of lines I]\label{local.genus.say}

Let $C\subset \p^n$ be a union of $m$ lines $L_i$.  We compute its
Hilbert polynomial in 2 ways.  Let $I\subset k[x_0,\dots, x_n]$
be the ideal of all homogeneous polynomials that vanish on $C$.
Then, for $d\gg 1$, $H_C(d)= md+1-p_a(C)$ is the dimension of the quotient
$$
W_C(d):=
\frac{(\mbox{degree $d$ homogeneous polynomials in $k[x_0,\dots, x_n] $})}
{ (\mbox{degree $d$ homogeneous polynomials that vanish on $C$})}.
$$
Let $W_{\p^n}(d)$ denote the  vector space of
degree $d$ homogeneous polynomials on $\p^n$. 
Let $\bar C:=\coprod L_i$ denote 
 the disjoint union of the lines $L_i$ and
 $\pi:\bar C\to C$ the natural map.

A degree $d$ homogeneous polynomial in $k[x_0,\dots, x_n] $
restricts to a degree $d$ homogeneous polynomial on each $L_i\cong \p^1$. 
This gives a restriction map
$$
\rest_d: W_{\p^n}(d)\into \tsum_{i=1}^m W_{L_i}(d)\cong \tsum_{i=1}^m W_{\p^1}(d)\cong k^{m(d+1)}
\eqno{(\ref{local.genus.say}.2)}
$$
that induces an injection 
$$
W_C(d)\to\tsum_{i=1}^m W_{L_i}(d).
\eqno{(\ref{local.genus.say}.3)}
$$
The linear terms of the Hilbert polynomials
 of the two sides of (\ref{local.genus.say}.3) are equal, hence we conclude that
$$
p_a(C)= \dim \bigl(\coker(\rest_d)\bigr)-m+1\qtq{for $d\gg 1$.}
\eqno{(\ref{local.genus.say}.4)}
$$
We aim to rewrite the cokernel of $\rest_d $ in terms of 
intersection points of the lines.
We start with the special case when there is only 1 intersection point.
The general formula will then be just a sum of such local terms.
\end{say}

\begin{say}[Local genus formula]\label{local.genus.say2}
Let $C_r^n\subset \a^n\subset \p^n$ be a union of $r$ lines $L_i$ through the origin.

A (parametrized, affine) line is given by $q:t\mapsto (a_1t,\dots, a_nt)$
and the corresponding restriction map is
$q^*: f(x_1,\dots, x_n)\mapsto  f(a_1t,\dots, a_nt)$.
Given $r$ different lines 
through the origin corresponds to  $r$ maps $q_i^*$.
A homogeneous polynomial of degree $d$ on $\p^n$ can be identified with a
polynomial of degree $\leq d$ on $\a^n$.
Thus the cokernel of $\rest_d$ is identified with
$$
\dim \coker\Bigl(k[x_1,\dots, x_n]_{(d)}
\stackrel{\oplus q_i^*}{\longrightarrow}\oplus_{i=1}^r k[t_i]_{(d)}\Bigr)
\eqno{(\ref{local.genus.say2}.1)}
$$
where the subscript $(d)$ denotes the subspace of polynomials of degree $\leq d$.

It is best to think of $k[x_1,\dots, x_n]_{(d)} $ as the
vector space of degree $\leq d$ Taylor polynomials on $0\in \a^n$ and
$k[t_i]_{(d)} $ as the
vector space of degree $\leq d$ Taylor polynomials on $0\in L_i$.

Fix a line $L_1$. For every other  line $L_i$
pick a linear form $\ell_i$ that vanishes on $L_i$ but not on $L_1$.
Set $g_1=\prod_{i>1}\ell_i$. Thus the image of $g_1$ under
 $\oplus q_i^* $ is zero in the summands $k[t_i]$ for $i>1$ 
and equals  $(\mbox{non-zero constant})\cdot t_1^{r-1}$
in $k[t_1]$. 
Therefore the cokernel of $\oplus q_i^* $
stabilizes for $d\geq r-1$.
This gives the local genus formula
$$
\delta(0\in C_r^n):=\dim \coker\bigl(k[x_1,\dots, x_n]_{(d)}
\stackrel{\oplus q_i^*}{\longrightarrow}\oplus_{i=1}^r k[t_i]_{(d)}\cong k^{r(d+1)}\bigr)
\eqno{(\ref{local.genus.say2}.2)}
$$
which holds for all $d\geq r-1$.

Since the $q_i^*$ preserve the degree, we can compute the cokernel
one degree at a time. 
In degree 0 there are just the constants
in $k[x_1,\dots, x_n] $ but $r$ copies of the constants
in the target in (\ref{local.genus.say2}.2). Thus
$$
\delta(0\in C_r^n)\geq r-1.
\eqno{(\ref{local.genus.say2}.3)}
$$
This leads to the weakest estimate  
(\ref{ci.genus.prop}.2).

In degree $j$ we have $\binom{j+n-1}{n-1}$ monomials of degree
$i$ in $k[x_1,\dots, x_n] $ and $r$ copies of $t_i^j$ in the target.
Therefore
$$
\delta(0\in C_r^n)\geq  \tsum_j \Bigl[r-\tbinom{j+n-1}{n-1}\Bigr]
\eqno{(\ref{local.genus.say2}.4)}
$$
where we sum over those $j\geq 0$ for which
the quantity in the brackets is positive.
(It is not hard to see that equality holds if the lines are in general position,
but this is not important for us.)
For $n=2$ this sum can be easily computed and we get that
$$
\delta(0\in C_r^2)=\tbinom{r}{2}.
\eqno{(\ref{local.genus.say2}.5)}
$$
(See \cite[Sec.IV.4.1]{shaf} for a different way of computing this.)
This leads to the strongest estimate  
(\ref{ci.genus.prop}.4).

If $n=3$ then there is no convenient closed form and the precise values
depend on the position of the lines. 
 For small vales of $r$ we get
$\delta(0\in C_2^3)\geq 1$, $\delta(0\in C_3^3)\geq 2$, $\delta(0\in C_4^3)\geq 4$ and 
$\delta(0\in C_5^3)\geq 6$.
It is easy to show that 
$$
\delta(0\in C_r^3)\geq  \tfrac1{\sqrt{2}} \bigl(r-1\bigr)^{3/2},
\eqno{(\ref{local.genus.say2}.6)}
$$
with equality holding only for $r=3$.

\end{say}

\begin{say}[Arithmetic genus of a union of lines II]\label{local.genus.say3}

Continuing the discussion of (\ref{local.genus.say}), 
pick any singular point $p\in C$. Let $C(p)\subset C$ denote the union of
 the lines  passing though $p$. As in (\ref{local.genus.say2}),
after choosing an affine chart and coordinates
we get  maps
between the spaces of Taylor polynomials
$$
\operatorname{Taylor}_d(p\in \a^n)
\stackrel{\oplus q_i^*}{\longrightarrow}  \oplus_i  
\operatorname{Taylor}_d(p\in L_i).
\eqno{(\ref{local.genus.say3}.1)}
$$
whose cokernel has dimension  $\delta(p\in C(p))$.
We can sum these over all singular points $\sing C$ of $C$ to get maps
$$
\operatorname{LocRest}_d: W_{\p^n}(d)\to 
 \bigoplus_{p\in \sing C} \oplus_i  \operatorname{Taylor}_d(p\in L_i).
\eqno{(\ref{local.genus.say3}.2)}
$$
Note that $ \operatorname{LocRest}_d$ factors through
$\rest_d$. Indeed, the map from
$W_{L_i}(d)$ to the right hand side of (\ref{local.genus.say3}.2) is obtained 
by starting with a
degree $d$ homogeneous polynomial $h$ on the line $L_i$ and for each
singular point $p\in L_i$ sending it to the degree $\leq d$ part of its
Taylor expansion at $p$. Each line contains at most $r-1$ singular points
thus these maps are surjective for $d> (r-1)^2$. This shows that
$$
\dim \coker \operatorname{LocRest}_d\geq \tsum_{p\in \sing C} \delta(p\in C(p)).
\eqno{(\ref{local.genus.say3}.3)}
$$
Combining (\ref{local.genus.say3}.3) with the local bounds
(\ref{local.genus.say2}.3--6) completes the  proof of 
Proposition \ref{ci.genus.prop}
for unions of lines once we prove  Proposition \ref{ag.prop.2}.

\end{say}

\begin{rem} For any (proper, reduced) algebraic curve $C$
there is a similar formula for the difference between
 the arithmetic genus of $C$
and the arithmetic genus of its normalization $\bar C$
in terms of local invariants computable from the singular points.
(These local terms are denoted by $\ell(\bar{\o}_x/\o_x)$ in   
\cite[Vol.1,p.262]{shaf}.) 

 If $(p\in C)$ is an analytically irreducible
curve singularity of multiplicity $r$ in $\c^3$ then 
$\delta(p\in C)\geq \rdown{r^2/4}$. Thus singularities with smooth branches
have the smallest genus for fixed multiplicity.
\end{rem}

\section{Arithmetic genus of subcurves}\label{sec.subcurves}

The  proof of Proposition \ref{ag.prop.2}
 uses basic sheaf cohomology theory. 
Everything we need is in  Sections III.1--5 of \cite{hartsh},
though the key statements are exercises.

First we use the cohomological
interpretation of the constant term of the Hilbert polynomial
as the holomorphic Euler characteristic.  This is a short argument.

\begin{lem} \label{hilb.char.lem}\cite[Exrc.III.5.2]{hartsh}
Let $I\subset k[x_0,\dots, x_n]$ be a homogeneous ideal such that
the corresponding scheme $C:=V(I)\subset \p^n$
 is 1-dimensional.
Then
\begin{enumerate}
\item  $h^0(C, \o_C)-h^1(C, \o_C)= H_C(0)$ and hence
\item $p_a(C)= h^1(C, \o_C)-h^0(C, \o_C)+1$. \qed
\end{enumerate}
\end{lem}

For complete intersection curves
we need  the following; this is a longer exercise.

\begin{lem} \label{hilb.char.lem2}\cite[Exrc.III.5.5]{hartsh}
For $i=1,\dots, n-1$ let $H_i\subset \p^n$ be 
a hypersurface of degree $a_i$.
Assume that the  intersection
$B:=H_1\cap \cdots\cap H_{n-1}$ is 1-dimensional.
Then
\begin{enumerate}
\item  $h^0(B, \o_B)=1$ and hence
\item $p_a(B)= h^1(B, \o_B)$. \qed
\end{enumerate}
\end{lem}

\begin{say}[Proof of Proposition \ref{ag.prop.2}]
We have a scheme theoretic intersection $B$ and
a reduced subcurve 
 $C\subset B$ which is defined by an
ideal sheaf $J_C\subset \o_B$. The exact sequence
$$
0\to J_C\to \o_B\to \o_C\to 0
$$
gives
$$
H^1(B, \o_B)\to H^1(C, \o_C)\to H^2(B, J_C)=0;
$$
the last vanishing holds since  $H^2$ is always zero on a curve;
cf.\ \cite[III.2.7]{hartsh}.
Thus  
$h^1(C, \o_C)\leq h^1(B, \o_B)$.

Since $C$ is reduced, $h^0(C,\o_C)$ equals 
the number of connected components of $C$.
Thus, by Lemma \ref{hilb.char.lem}.2 and Lemma \ref{hilb.char.lem2}.2,
$$
p_a(C)=h^1(C, \o_C)-h^0(C, \o_C)+1\leq h^1(C, \o_C)\leq 
 h^1(B, \o_B)=p_a(B).\qed
$$
\end{say}

\begin{rem} I tried to find a more elementary 
proof of Proposition \ref{ag.prop.2} but so far I have been unsuccessful.
There is a vast classical literature on curves in $\p^3$, but most of
it studies smooth or only mildly singular curves. 

Let $X$ be a normal, projective variety of dimension $n$
and $H_1,\dots, H_{n-1}$ hyperplane sections such that
$B:=H_1\cap \cdots \cap H_{n-1}$ is 1-dimensional. There is a
formula similar to (\ref{hilb.say}.5) that computes the genus
of $B$ if $B$ is smooth. However, when $B$ is singular and
$C\subset B$ is a reduced subcurve, it can happen that
the arithmetic genus of $C$ is bigger than the
arithmetic genus of $B$. Thus Proposition \ref{ag.prop.2}
is a special property of $\p^n$. However, in all the examples that I computed,
the arithmetic genus of $C$ is not much bigger than the
arithmetic genus of $B$.

Chasing through the proofs of (\ref{hilb.char.lem2}),
 the key property seems to be that Kodaira's vanishing theorem
holds for  $\p^n$.
\end{rem}

\section{Ruled surfaces}
\label{ruled.sec}

The referee pointed out that information about ruled surfaces is hard to
extract from the current literature, so here I summarize the pertinent facts
with proofs.

We are interested in the geometry of ruled surfaces, thus in this section
we work over an algebraically closed field  $K$, though almost
everything works over any infinite field.

\begin{defn}
A {\it smooth minimal ruled surface} is a projective surface $M$
with a morphism $g:M\to C$ to a smooth curve all of  whose fibers,
also called {\it rulings,}
are  (isomorphic to) lines.


 A {\it ruled surface} is a projective surface $S\subset \p^n$
that is the image of a  smooth ruled surface $M$
by a  morphism $\pi:M\to S$ that sends the rulings to lines.
We call  
$$
C\stackrel{g}{\leftarrow}M \stackrel{\pi}{\to} S
$$
a {\it presentation} of $S$. 
If $\pi$ is birational, we call it a
{\it birational presentation.}
Any surface in $\p^n$ can be birationally projected to $\p^3$,
so we focus on surfaces in $\p^3$. 

We will show that every ruled surface has a  
 birational presentation, and, with two exceptions, the
 birational presentation is unique. Thus birationality is frequently
part of the definition. (Note  that the literature is inconsistent.
Sometimes a ruled surface means a smooth minimal ruled surface,
a ruled surface as above or any  surface that is birational to
a ruled surface.)

\end{defn}

\begin{prop} \label{lines.on.surf.prop}
Let $S\subset \p^3$ be an irreducible  surface of degree $d$. Then
\begin{enumerate}
\item either  $S$ contains at most $d^4$ lines
\item or $S$ is ruled.
\end{enumerate}
\end{prop}

Proof. First we use affine coordinates. A typical
 line on $S$ can be given parametrically as
$t\mapsto  (a_1t+b_1, a_2t+b_2, t)$.
If $f(x,y,z)=0$ is an affine equation of $S$, such a line is
contained in $S$ iff
$$
f\bigl(a_1t+b_1, a_2t+b_2, t\bigr)\equiv 0.
$$
Expanding by the powers of $t$, we get a system of  $d+1$
equations of degree $\leq d$ in the variables $a_1, b_1, a_2, b_2$.
By B\'ezout, the system either has at most $d^4$ solutions
(leading to the first case) or the solution set contains an algebraic curve
$C\subset \a^4$ (with $a_1, b_1, a_2, b_2$ as coordinates). In this case
$$
\pi: C\times \a^1\to S\qtq{given by} (c,t)\mapsto (a_1t+b_1, a_2t+b_2, t)
$$
is a rational map from an (affine)
 ruled surface to $S$. There could be several such curves $C$
and the resulting map $\pi$ need not be birational, but we do get
at least 1 rational presentation of $S$.

A few details need to be ironed out. In general, $C$ is neither smooth nor
projective. Thus one should work with the Grassmannian parametrizing
all  lines in $\p^3$ ; see \cite[Vol.1,p.42]{shaf}. 
Then we have to normalize $C$
to get a smooth ruled surface mapping onto $S$.  
\qed

\begin{say}[Special ruled surfaces] \label{ruled.specv.say}
There are 3 types of ruled surfaces that are exceptional for many
of the  results. These are
\begin{enumerate}
\item planes,
\item smooth quadrics and
\item cones.
\end{enumerate}
The plane has infinitely many birational presentations given by the family
of all lines passing through a given point. Correspondingly, the plane
can be viewed as a cone in infinitely many ways.
A smooth quadric has 2 birational presentations.
A cone (that is not a plane) has a unique  birational presentation but all the lines
pass through the unique vertex. 
Every other ruled surface will be called {\it non-special.}

Fix a presentation  $C\stackrel{g}{\leftarrow}M \stackrel{\pi}{\to} S$ and
 let $Z\subset S$ be any subset.
Then $g\bigl(\pi^{-1}(Z)\bigr)\subset C$ is the set of rulings that meet
 $Z$ in at least 1 point.

If $B=Z$ is an irreducible
 curve then $\pi^{-1}(B)\subset M$ is also a curve hence
\begin{enumerate}\setcounter{enumi}{3}
\item  either
$g\bigl(\pi^{-1}(B)\bigr)\subset C$ consist of finitely many points,
$B$ is a ruling and  only finitely many rulings of the given presentation
 intersect $B$, 
\item or $g\bigl(\pi^{-1}(B)\bigr)= C$, hence every ruling
intersects $B$. 
\end{enumerate}
A line $L\subset S$ is called {\it special}
(for the given presentation) if it intersects every ruling.

A point $p\in S$ is  called {\it special} if it is either singular or it 
lies on a special line. (We will show that a non-special ruled surface
has non-special points.)
\end{say}

\begin{prop} \label{ruled.surf.long.prop}
Let $S\subset \p^3$ be a non-special ruled surface of degree $d$. Then
\begin{enumerate}
\item there are at most $d$ lines through any point of $S$,
\item there is exactly one line  through a  non-special point of $S$,
\item a non-special ruled surface admits a unique birational presentation,
\item there are at most 2 special lines and
\item there are at most $d-2$
  non-special lines intersecting a non-special line.
\end{enumerate}
\end{prop}

Proof. Assume that $p\in S$ is a point with infinitely any lines through it.
Choose affine coordinates such that $p=(0,0,0)$ and $S$ has equation
$f(x,y,z)=\sum_i f_i(x,y,z)$ where  $f_i$ is homogeneous of degree $i$.
A parametrized line $t\mapsto (at, bt,ct)$ lies on $S$ iff
$f( at, bt,ct)$ is identically 0. This holds iff
$f_i(a,b,c)=0$ for every $i$. By B\'ezout, there are either
finitely many (in fact $\leq d(d-1)$) solutions or the $f_i$ have a common 
(homogeneous) factor
 $h(x,y,z)$. Then $h$ divides $f$ hence  the cone
$(h=0)$ is an irreducible component of $S$. 
This is a contradiction since $S$ is irreducible and not a cone.
Thus $\pi:M\to S$ is everywhere finite--to--one.

Next we claim that any 2 special lines  $L_1, L_2$ are disjoint.
If not then they span a plane $P$. As we noted, only 
finitely many rulings pass through the point $L_1\cap L_2$,
hence every other ruling meets $P$ in 2 points. Thus every other ruling
is contained in $P$ hence $S=P$.

Assume next that $S$ contains 3 special lines  $L_1, L_2, L_3$.
For a quadric it is 3 conditions to contain a line, hence there is
a quadric $Q$ that contains all 3 lines. 
Thus every  ruling  meets $Q$ is at least 3 points hence is contained in
$Q$. Thus $S=Q$, proving (4).

Let $p\in S$ be a smooth point and $p\in B\subset S$
a line. Since $p$ is smooth, $B$ is locally defined by
1 equation at $p$ \cite[Vol.1,p.108]{shaf}, 
thus $\pi^{-1}(B)$ is locally defined by
1 equation at $\pi^{-1}(p)$. Thus  $\pi^{-1}(B)$ is 1-dimensional at
$\pi^{-1}(p)$ \cite[Vol.1,p.71]{shaf}.
 Thus either $B$ is a ruling passing through $p$
or $B$ is special. Thus there is exactly one line  through a 
 non-special point of $S$,
proving (2). Hence in the construction of Proposition \ref{lines.on.surf.prop}
the curve $C$ is unique and the resulting $\pi:M\to S$
is the unique birational ruling of $S$, proving (3). 
(Strictly speaking, we have only proved that $\pi:M\to S$ is injective
on a dense open subset. This implies birationality in characteristic 0.
In positive characteristic  we still need to exclude purely inseparable maps.
Since this has no bearing on curve counts, we do not pursue this issue.)

In order to get precise bounds on the number of lines, we
use intersection theory (\ref{int.thm.say}) on the smooth surface $M$
for the family  $\{H_{\lambda}\}$ of pull-backs of plane sections
$S\cap P_{\lambda}\subset S$ of $S$.

First choose planes $P_1, P_2\subset \p^3$ such that the line $P_1\cap P_2$
meets $S$ in $d$ distinct smooth points. Then $H_1$ and $ H_2$
meet at the preimages of these points and $m_p(H_1,H_2)=1$ at each of them.
Thus $(H_1\cdot H_2)=d$. 

Next let $p\in S$ be any point and  choose planes $P_1, P_2$ such that the 
line $P_1\cap P_2$
meets $S$ at $p$ but is not contained in $S$. 
As we noted in (\ref{ruled.specv.say}),
the rulings passing through $p$ correspond to the set
$g\bigl(\pi^{-1}(p)\bigr)$, hence its cardinality is at most
$|\pi^{-1}(p)|=|H_1\cap H_2|$. By (\ref{ruled.surf.long.prop}.6),
$|H_1\cap H_2|\leq (H_1\cdot H_2)=d$. 
Thus there are at most $d$ rulings passing through $p$.
We can arrange that  $P_1\cap P_2$ meets $S$ in at least one more point;
this shows that there are at most $d-1$ rulings passing through $p$
and at most 1 special line, proving (1).

Finally let $L\subset S$ be a ruling and choose $P_1, P_2$ such that  
$L=P_1\cap P_2$.
The corresponding $H_1, H_2$ are reducible curves of the form
$H_i=B_i+C_i$ where  $B_i=\sum_j a_{ij}F_j$, $a_{ij}>0$ and 
$F_j\subset M$ are the rulings   such that 
$L=\pi(F_j)$. The other rulings
intersecting $L$ correspond to the points $C_1\cap C_2$. 
As before, we compute the intersection number
$$
\begin{array}{rcl}
|C_1\cap C_2|&\leq &(C_1\cdot C_2)\\
&=&\bigl((H_1-B_1)\cdot (H_1-B_2)\bigr)\\
&=&
 (H_1\cdot H_2)- (H_1\cdot B_2)- (B_1\cdot H_2)+(B_1\cdot B_2)\\
&\leq & d-2. 
\end{array}
$$
This proves (5). \qed

\begin{say}[Intersecting curves on smooth surfaces]\label{int.thm.say}
For proofs see  \cite[Sec.IV.1]{shaf} or \cite[V.1]{hartsh}.

Let $X$ be a smooth, projective surface. 
 Given two curves $A,B\subset X$, there is an intersection
number $(A\cdot B)$ attached to them. This number is
symmetric, bilinear and unchanged if we vary the curves in
families. Furthermore, if $A\cap B$ is finite
then 
$$
(A\cdot B)=\sum_{p\in A\cap B}m_p(A,B)
\eqno{(\ref{int.thm.say}.1)}
$$
where each $m_p(A,B)$ is a positive integer. 
Furthermore  $m_p(A,B)=1$ iff $A,B$ are both smooth at $p$ and 
are not tangent there.
\end{say}

\section{Sketch of the proof of the Monge--Salmon--Cayley theorem}
\label{sketch.sec}

\begin{say}[Salmon's flecnodal equation] \label{salmon.surf.say}
Let us start with 3 homogeneous forms in 3 variables
$$
\tsum_{1\leq i\leq 3} a_i x_i, \quad \tsum_{1\leq i\leq j\leq 3} b_{ij} x_ix_j, 
\quad\tsum_{1\leq i\leq j\leq k\leq 3} c_{ijk} x_ix_ix_k.
\eqno{(\ref{salmon.surf.say}.1)}
$$
We want to understand when they have a common zero.
We eliminate $x_3$ from the linear equation and substitute into the others
to get 2 homogeneous forms in 2 variables
$$
\tsum_{1\leq i\leq j\leq 2} B_{ij} x_ix_j, \quad\tsum_{1\leq i\leq j\leq k\leq 2} C_{i j k} x_ix_jx_k.
\eqno{(\ref{salmon.surf.say}.2)}
$$
They have a common zero iff their discriminant vanishes.
After clearing the denominator (which is a power of $a_3$)
this gives an equation in the original variables  $a_i,b_{ij}, c_{ijk}$.
After a short argument about the $a$-variables we get the following.
\medskip

{\it Claim} \ref{salmon.surf.say}.3. There is a polynomial
$F(\ , \ , \ )$ such that $F(a_i,  b_{ij}, c_{ijk})=0$
iff the 3 forms in (\ref{salmon.surf.say}.1) have a common (nontrivial)
 zero.
Furthermore, $F$ has multidegree  $(6,3,2)$. \qed
\medskip

Consider now a surface  $S\subset \c^3$ given by an equation
$f(x_1, x_2, x_3)=0$. Fix a point $p=(p_1, p_2, p_3)\in S$
and write  the Taylor expansion of $f$ around $p$
as
$$
f=\tsum_{i=0}^d f_i(x_1-p_1, x_2-p_2, x_3-p_3)
\eqno{(\ref{salmon.surf.say}.4)}
$$
where $f_i$ is homogeneous of degree $i$. A parametric line
$$
t\mapsto (p_1+m_1t, p_2+m_2t, p_3+m_3t)
$$ is a triple tangent
iff
$$
f_1(m_1,m_2,m_3)=f_2(m_1,m_2,m_3)=f_3(m_1,m_2,m_3)=0.
\eqno{(\ref{salmon.surf.say}.5)}
$$
By (\ref{salmon.surf.say}.3) this translates into 
an equation $F(a_i,  b_{ij}, c_{ijk})=0$ in the
coefficients of the $f_i$, which are in turn given by the
$i$th partial derivatives of $f$.

Putting all together we get a polynomial
$$
\operatorname{Flec}_f(x_1, x_2, x_3):=F\Bigl( \frac{\partial f}{\partial x_i}, 
\frac{\partial^2 f}{\partial x_i\partial x_j}, 
 \frac{\partial^3 f}{\partial x_i\partial x_j\partial x_k}\Bigr)
\eqno{(\ref{salmon.surf.say}.6)}
$$
such that
$$
f(x_1, x_2, x_3)=\operatorname{Flec}_f(x_1, x_2, x_3)=0
\eqno{(\ref{salmon.surf.say}.7)}
$$
defines the set of points of $S$ where there is a triple tangent line.
Furthermore,  $\operatorname{Flec}_f $
 has degree  $\leq 6(d-1)+3(d-2)+2(d-3)=11d-18$ in $x,y,z$. 

Note that the coefficients of the different $f_i$
are not independent, thus one could end up with a lower
degree polynomial. Salmon claims that in fact one gets a
polynomial of degree $11d-24$. 
I have not checked this part; in our applications we have 
used only that the degree is $\leq 11d$.

Note that when $\deg f=3$, the Salmon bound is
$11\cdot 3-24=9$. A smooth cubic surface $S$ contains 27 lines
and their union is the complete intersection of $S$ with
a surface $T$ of degree $9$. So, in this case, the
Salmon bound is sharp.
\end{say}

If a line is contained in $S$, then it is  triply tangent everywhere,
thus  $\operatorname{Flec}_f $ vanishes on every line contained in $S$.
This is useful only if  $\operatorname{Flec}_f $ does not vanish
identically on $S$. That is, we need to understand surfaces
where every point has a triple tangent line. 
Monge proved that these are exactly the ruled surfaces.
Monge writes a surface locally as a graph, thus from now on
we work with holomorphic functions (over $\c$)
or with $C^3$-functions (over $\r$). 

\begin{say}[Monge's theorem] \label{monge.say}
Consider a graph $S:=\bigl(z=f(x,y)\bigr)\subset \c^3$.
Fix a point $(x_0, y_0, z_0)$. The line
$$
(x_0+t, y_0+mt, z_0+nt)
\eqno{(\ref{monge.say}.1)}
$$
is a double tangent line of $S$ iff
$n=f_x(x_0, y_0)+f_y(x_0, y_0)m$ and 
$$
f_{xx}(x_0, y_0)+2f_{xy}(x_0, y_0)m+f_{yy}(x_0, y_0)m^2=0.
\eqno{(\ref{monge.say}.2)}
 $$
The double tangent lines are also called {\it asymptotic directions.}
By working on a smaller open set,
we may assume that the Hessian of $f$ has constant rank and is not 
identically 0. Thus the asymptotic directions define 2 vector fields
on $S$. (Only 1  vector field if the rank is always 1.)
Integrating these vector fields we get the
{\it asymptotic curves} of the surface $S$.

The line (\ref{monge.say}.1) is a triple tangent
if, in addition
$$
f_{xxx}(x_0, y_0)+3f_{xxy}(x_0, y_0)m+3f_{xyy}(x_0, y_0)m^2+f_{yyy}(x_0, y_0)m^3=0.
\eqno{(\ref{monge.say}.3)}
 $$
Thus the graph has a triple tangent iff the equations
(\ref{monge.say}.2--3) have a  common solution.

\medskip

{\it Claim} \ref{monge.say}.4. An asymptotic curve is a straight line
iff all the corresponding asymptotic directions are triple tangents.
\medskip

Proof. 
Assume that we have $u=u(t)$ defined by
$a(t)+2b(t)u+ c(t)u^2=0$.
By implicit differentiation,  $u(t)$ is constant iff
$a_t+2b_tu+ c_tu^2\equiv 0$.
 
Assume next that  $u=u(x,y)$ is defined by
$$
a(x,y)+2b(x,y)u+ c(x,y)u^2=0
$$
and we work along a path  $\bigl(x(t), y(t)\bigr)$.
Then the condition becomes
$$
a_xx'+\bigl(a_yy'+2b_xux'\bigr)+
\bigl(2b_yuy'+c_xu^2x'\bigr)+ c_yu^2y'\equiv 0.
$$
In our case  $a=f_{xx}, b=f_{xy}, c=f_{yy}$  and $u=y'/x'$
along the asymptotic curve. 
Substituting $y'=ux'$ and dividing by $x'$
we get the condition
$$
f_{xxx}+3f_{xxy}u+3f_{xyy}u^2+f_{yyy}u^3=0,
$$
which is the same as (\ref{monge.say}.3).\qed
\end{say}

See  \cite[2.10]{MR739791} or 
\cite{tao-blog-monge} for other variants of this argument.

\def\cprime{$'$} \def\cprime{$'$} \def\cprime{$'$} \def\cprime{$'$}
  \def\cprime{$'$} \def\cprime{$'$} \def\dbar{\leavevmode\hbox to
  0pt{\hskip.2ex \accent"16\hss}d} \def\cprime{$'$} \def\cprime{$'$}
  \def\polhk#1{\setbox0=\hbox{#1}{\ooalign{\hidewidth
  \lower1.5ex\hbox{`}\hidewidth\crcr\unhbox0}}} \def\cprime{$'$}
  \def\cprime{$'$} \def\cprime{$'$} \def\cprime{$'$}
  \def\polhk#1{\setbox0=\hbox{#1}{\ooalign{\hidewidth
  \lower1.5ex\hbox{`}\hidewidth\crcr\unhbox0}}} \def\cdprime{$''$}
  \def\cprime{$'$} \def\cprime{$'$} \def\cprime{$'$} \def\cprime{$'$}
\providecommand{\bysame}{\leavevmode\hbox to3em{\hrulefill}\thinspace}
\providecommand{\MR}{\relax\ifhmode\unskip\space\fi MR }
\providecommand{\MRhref}[2]{%
  \href{http://www.ams.org/mathscinet-getitem?mr=#1}{#2}
}
\providecommand{\href}[2]{#2}

\vskip1cm

\noindent Princeton University, Princeton NJ 08544-1000

{\begin{verbatim}kollar@math.princeton.edu\end{verbatim}}

\end{document}